\documentclass[9pt, letterpaper]{article}
\usepackage[latin1]{inputenc}
\usepackage[T1]{fontenc}
\usepackage[english]{babel}
\usepackage{amsmath}
\usepackage{amssymb}
\usepackage{hyperref}

\usepackage[margin = 1in]{geometry}

\allowdisplaybreaks

\newtheorem{Lemma}{Lemma}
\newtheorem{Theorem}{Theorem}
\newtheorem{Remark}{Remark}
\newenvironment{Proof}{\textbf{Proof:}}{\hfill $\Box$}

\begin{document}

\title{Nonlinear damped Timoshenko systems with second sound --- global existence and exponential stability}
\author{Salim A. Messaoudi\thanks{%
	Mathematical Sciences Department, KFUPM, Dhahran 31261, Saudi Arabia\newline
	E-mail: messaoud@kfupm.edu.sa}, \thinspace Michael Pokojovy\thanks{%
	Fachbereich Mathematik und Statistik, Universit\"{a}t Konstanz, 78457 Konstanz,
	Germany\newline
	E-mail: Michael.Pokojovy@uni-konstanz.de}, \thinspace Belkacem Said-Houari%
	\thanks{%
	Universit\'{e} Badji Mokhtar, Laboratoire de Math\'{e}matiques Appliqu\'{e}es, B.P. 12
	Annaba 23000, Algerie\newline
	E-mail: saidhouarib@yahoo.fr}}
	\date{March 2008}
	\maketitle

\begin{abstract}
	In this paper, we consider nonlinear thermoelastic systems of Timoshenko
	type in a one-dimensional bounded domain. The system has two dissipative
	mechanisms being present in the equation for transverse displacement and
	rotation angle --- a frictional damping and a dissipation through hyperbolic
	heat conduction modelled by Cattaneo's law, respectively. The global
	existence of small, smooth solutions and the exponential stability in linear
	and nonlinear cases are established.
\end{abstract}

\vspace{0.3cm}

\noindent AMS-Classification: 35B37, 35L55, 74D05, 93D15, 93D20 \newline

\section{Introduction}
	In \cite{Ti1921}, a simple model describing the transverse vibration of a
	beam was developed. This is given by a system of two coupled hyperbolic
	equations of the form 
	\begin{eqnarray}
		\rho u_{tt} &=& (K(u_x - \varphi))_x \quad \text{ in } (0, \infty) \times (0, L),  \label{TIMOSHENKO_LIN_1} \\
		I_\rho \varphi_{tt} &=& (EI \varphi_x)_x + K(u_x - \varphi) \quad \text{ in } (0, \infty) \times (0, L),  \nonumber
	\end{eqnarray}
	where $t$ denotes the time variable and $x$ the space variable along a beam
	of length $L$ in its equilibrium configuration. The unknown functions $u$
	and $\varphi$ depending on $(t, x) \in (0, \infty) \times (0, L)$ model the
	transverse displacement of the beam and the rotation angle of its filament,
	respectively. The coefficients $\rho$, $I_\rho$, $E$, $I$ and $K$ represent
	the density (i.e. the mass per unit length), the polar momentum of inertia
	of a cross section, Young's modulus of elasticity, the momentum of inertia
	of a cross section, and the shear modulus, respectively.
	
	Kim and Renardy considered (\ref{TIMOSHENKO_LIN_1}) in \cite{KiRe1987} together with two boundary controls of the form 
	\begin{eqnarray}
		K\varphi(t, L) - Ku_{x}(t,L) &=& \alpha u_{t}(t, L) \quad \text{ in }(0,\infty), \nonumber \\
		EI\varphi_{x}(t, L) &=& -\beta \varphi_{t}(t, L) \quad \text{ in }(0,\infty) \nonumber
	\end{eqnarray}
	and used the multiplier techniques to establish an exponential decay result
	for the natural energy of (\ref{TIMOSHENKO_LIN_1}). They also provided some
	numerical estimates to the eigenvalues of the operator associated with the
	system (\ref{TIMOSHENKO_LIN_1}). An analogous result was also established by
	Feng \textit{et al.} in \cite{FeShiZha1998}, where a stabilization of
	vibrations in a Timoshenko system was studied. Rapos \textit{et al.} studied
	in \cite{RaFeSaCa2005} the following system 
	\begin{align}
		& \rho_{1} u_{tt} - K(u_{x} - \varphi)_{x} + u_{t} = 0 \quad \text{ in } (0,\infty) \times (0, L),  \nonumber \\
		& \rho_{2} - b\varphi_{xx} + K(u_{x} - \varphi) + \varphi_{t} = 0 \quad \text{ in } (0,\infty) \times (0,L),  \label{TIMOSHENKO_DAMPED_TWICE} \\
		& u(t, 0) = u(t, L) = \varphi(t, 0) = \varphi(t, L) = 0 \quad \text{ in }(0, \infty) \nonumber
	\end{align}
	and proved that the energy associated with (\ref{TIMOSHENKO_DAMPED_TWICE})
	decays exponentially. This result is similar to that one by Taylor \cite{Ta2000}, 
	but as they mentioned, the originality of their work lies in the
	method based on the semigroup theory developed by Liu and Zheng \cite{LiuZhe1999}.
	
	Soufyane and Wehbe considered in \cite{SouWeh2003} the system 
	\begin{align}
		&\rho u_{tt} = (K(u_x - \varphi))_x \quad \text{ in } (0, \infty) \times (0,L),  \nonumber \\
		&I_\rho \varphi_{tt} = (EI\varphi_x)_x + K(u_x - \varphi) - b\varphi_t \quad 
		\text{ in } (0, \infty) \times (0, L),  \label{TIMOSHENKO_DAMPLED_ONCE} \\
		&u(t, 0) = u(t, L) = \varphi(t, 0) = \varphi(t, L) = 0 \quad \text{ in } (0, \infty),  \nonumber
	\end{align}
	where $b$ is a positive continuous function satisfying 
	\begin{equation}
		b(x) \geq b_0 > 0 \quad \text{ in } [a_0, a_1] \subset [0, L].  \nonumber
	\end{equation}
	In fact, they proved that the uniform stability of (\ref
	{TIMOSHENKO_DAMPLED_ONCE}) holds if and only if the wave speeds are equal, i.e. 
	\begin{equation}
		\frac{K}{\rho} = \frac{EI}{I_\rho},  \nonumber
	\end{equation}
	otherwise, only the asympotic stability has been proved. This result
	improves previous ones by Soufyane \cite{Sou1999} and Shi and Feng \cite{ShiFe2002} 
	who proved an exponential decay of the solution of (\ref{TIMOSHENKO_LIN_1}) 
	together with two locally distributed feedbacks.
	
	Recently, Rivera and Racke \cite{MuRa2007} obtained a similar result in a
	work where the damping function $b = b(x)$ is allowed to change its sign.
	Also, Rivera and Racke \cite{MuRa2003} treated a nonlinear Timoshenko-type
	system of the form 
	\begin{align}
		&\rho_1 \varphi_{tt} - \sigma_1(\varphi_x, \psi)_x = 0,  \nonumber \\
		&\rho_2 \psi_{tt} - \chi(\psi_x)_x + \sigma_2(\varphi_x, \psi) + d\psi_t = 0 \nonumber
	\end{align}
	in a one-dimensional bounded domain. The dissipation is produced here
	through a frictional damping which is only present in the equation for the
	rotation angle. The authors gave an alternative proof for a necessary and
	sufficient condition for exponential stability in the linear case and then
	proved a polynomial stability in general. Moreover, they investigated the
	global existence of small smooth solutions and exponential stability in the
	nonlinear case.
	
	Xu and Yung \cite{XuYu2003} studied a system of Timoshenko beams with
	pointwise feedback controls, looked for the information about the
	eigenvalues and eigenfunctions of the system, and used this information to
	examine the stability of the system.
	
	Ammar-Khodja \textit{et al.} \cite{AmBeRiRa2003} considered a linear
	Timoshenko-type system with a memory term of the form 
	\begin{align}
		& \rho_{1} \varphi_{tt} - K(\varphi_{x} + \psi)_{x} = 0, \label{TIMOSHENKO_MEMORY_TERM} \\
		& \rho_{2} \psi_{tt} - b \psi_{xx} + \int_{0}^{t} g(t - s) \psi_{xx}(s)ds + K(\varphi_{x} + \psi) = 0 \nonumber
	\end{align}
	in $(0,\infty) \times (0,L)$, together with homogeneous boundary conditions.
	They applied the multiplier techniques and proved that the system is
	uniformly stable if and only if the wave speeds are equal, 
	i.e. $\frac{K}{\rho_{1}} = \frac{b}{\rho_{2}}$, and $g$ decays uniformly. 
	Precisely, they proved an exponential decay if $g$ decays exponentially and polynomial decay
	if $g$ decays polynomially. They also required some technical conditions on
	both $g^{\prime}$ and $g^{\prime \prime}$ to obtain their result. The
	feedback of memory type has also been studied by Santos \cite{Sa2002}. He
	considered a Timoshenko system and showed that the presense of two feedbacks
	of memory type at a subset of the bounary stabilizes the system uniformly.
	He also obtained the energy decay rate which is exactly the decay rate of
	the relaxation functions.
	
	Shi and Feng \cite{ShiFe2001} investigated a nonuniform Timoshenko beam and
	showed that the vibration of the beam decays exponentially under some
	locally distributed controls. To achieve their goal, the authors used the
	frequency multiplier method.
	
	For Timoshenko systems of classical thermoelasticity, Rivera and Racke \cite
	{MuRa2007} considered, in $(0,\infty )\times (0,L)$, the following system 
	\begin{align}
		& \rho_{1} \varphi_{tt} - \sigma (\varphi_{x}, \psi_{x})_{x} = 0, \nonumber \\
		& \rho_{1} \psi_{tt} - b \psi_{xx} + k(\varphi_{x} + \psi) + \gamma \theta_{x} = 0, \label{TIMOSHENKO_NONLINEAR_FOURIER} \\
		& \rho_{3} \theta_{t} - \kappa \theta_{xx} + \gamma \psi_{tx} = 0,  \nonumber
	\end{align}
	where the functions $\varphi$, $\psi$, and $\theta$ depend on $(t, x)$ and
	model the transverse displacement of the beam, the rotation angle of the
	filament, and the temperature difference, respectively. Under appropriate
	conditions on $\sigma$, $\rho_{i}$, $b$, $k$, $\gamma $ they proved
	several exponential decay results for the linearized system and
	non-exponential stability result for the case of different wave speeds.
	
	In the above system, the heat flux is given by the Fourier's law. As a
	result, we obtain a physical discrepancy of infinite heat propagation speed.
	That is, any thermal disturbance at a single point has an instantaneous
	effect everywhere in the medium. Experiments showed that heat conduction in
	some dielectric crystals at low temperatures is free of this paradox.
	Moreover, the disturbances being almost entirely thermal, propagate at a
	finite speed. This phenomenon in dielectric crystals is called second sound.
	
	To overcome this physical paradox, many theories have been developed. One of
	which suggests that we should replace the Fourier's law 
	\begin{equation}
		q + \kappa \theta_x = 0  \nonumber
	\end{equation}
	by so called Cattaneo's law 
	\begin{equation}
		\tau q_t + q + \kappa \theta_x = 0.  \nonumber
	\end{equation}
	
	Few results concerning existence, blow-up, and asymptotic behavior of smooth
	as well as weak solutions in thermoelasticity with second sound have been
	established over the past two decades. Tarabek \cite{Ta1992} treated
	problems related to the following one-dimensional system 
	\begin{eqnarray}
		u_{tt} - a(u_{x}, \theta, q) u_{xx} + b(u_{x}, \theta, q) \theta_{x} &=& \alpha_{1}(u_{x}, \theta) qq_{x},  \nonumber \\
		\theta_{t} + g(u_{x}, \theta, q) q_{x} + d(u_{x}, \theta, q) u_{tx} &=& \alpha_{2}(u_{x}, \theta) qq_{t}, 
		\label{THERMOELASTICITY_CATTANEO_NONLINEAR_1D} \\
		\tau(u_{x}, \theta) q_{t} + q +k(u_{x}, \theta) \theta_{x} &=& 0  \nonumber
	\end{eqnarray}
	in both bounded and unbounded situations and established global existence
	results for small initial data. He also showed that these ``classical''
	solutions tend to equilibrium as $t$ tends to infinity. However, no decay
	rate has been discussed. Racke \cite{Ra2002} discussed lately (\ref{THERMOELASTICITY_CATTANEO_NONLINEAR_1D}) 
	and established exponential decay results for several linear and nonlinear initial boundary value problems. In
	particular, he studied the system (\ref{THERMOELASTICITY_CATTANEO_NONLINEAR_1D}) for a rigidly clamped medium with
	the temperature held constant on the boundary, i.e. 
	\begin{equation}
		u(t, 0) = u(t, 1) = \theta(t, 0) = \theta(t, 1) = \bar{\theta} \quad \text{ in } (0, \infty),  \nonumber
	\end{equation}
	and showed for sufficiently small initial data and $\alpha_{1} = \alpha_{2} = 0$ 
	that the classical solution decays exponentially to an equilibrium state. 
	Messaoudi and Said-Houari \cite{MeSa2005} extended the decay result of \cite{Ra2002} for $\alpha_{1}$ and $\alpha_{2}$ that are not
	necessarily zero.
	
	Concerning the multi-dimensional case ($n = 2,3$), Racke \cite{Ra2003}
	established an existence result for the following $n$-dimensional problem 
	\begin{align}
		& u_{tt} - \mu \Delta u - (\mu + \lambda) \nabla \textrm{\thinspace div \thinspace} u 
		+ \beta \nabla \theta = 0, \quad (t,x) \in (0, \infty) \times \Omega, \nonumber \\
		& \theta_{t} + \gamma \textrm{\thinspace div \thinspace} q 
		+ \delta \textrm{\thinspace div\thinspace} u_{t} = 0, \quad (t, x) \in (0, \infty) \times \Omega, \nonumber \\
		& \tau q_{t} + q + \kappa \nabla \theta = 0, \quad (t, x) \in (0, \infty) \times \Omega, \label{THERMOELASTICITY_CATTANEO_3D} \\
		& u(0, x) = u_{0}(x), \, u_{t}(0, x) = u_{1}(x), \, \theta(0, x) = \theta_{0}(x), \,q(0,x) = q_{0}(x), \quad x \in \Omega \nonumber \\
		& u(t, x) = \theta(t, x) = 0, \quad (t, x) \in (0, \infty) \times \partial \Omega, \nonumber
	\end{align}
	where $\Omega$ is a bounded domain of $\mathbb{R}^{n}$ with a smooth boundary $\partial \Omega $. 
	$u = u(t, x) \in \mathbb{R}^{n}$ is the displacement vector, $\theta = \theta(t, x)$ is the temperature difference, 
	$q = q(t, x) \in \mathbb{R}^{n}$ is the heat flux, and $\mu$, $\lambda$, $\beta$, $\gamma$, $\delta$, $\tau$, $\kappa$ are positive constants, 
	where $\mu$, $\alpha$ are Lam\'{e} moduli and $\tau$ is the relaxation time being a small parameter compared to the others. 
	In particular, if $\tau = 0$, the system (\ref{THERMOELASTICITY_CATTANEO_3D}) reduces to the system of thermoelasticity, 
	in which the heat flux is given by Fourier's law instead of Cattaneo's law. He also proved, under condition 
	$\nabla \times \nabla u = \nabla \times \nabla q = 0$, an exponential decay result for (\ref{THERMOELASTICITY_CATTANEO_3D}). 
	This result is easily extended to the radially symetric solutions, as they satisfy the above condition.
	
	Messaoudi \cite{Me2002} investigated the following problem 
	\begin{align}
		& u_{tt} - \mu \Delta u - (\mu + \lambda) \nabla \textrm{\thinspace div\thinspace} u
		+ \beta \nabla \theta = |u|^{p-2} u, \quad (t, x) \in (0, \infty) \times \Omega, \nonumber \\
		& \theta_{t} + \gamma \textrm{\thinspace div\thinspace } q + \delta \textrm{\thinspace div\thinspace } u_{t}
		= 0, \quad (t, x) \in (0, \infty) \times \Omega, \nonumber \\
		& \tau q_{t} + q + \kappa \nabla \theta = 0, \quad (t, x) \in (0, \infty) \times \Omega ,  \label{THERMOELATICITY_CATTANEO_POLYNOMIAL_3D} \\
		& u(0, x) = u_{0}(x), \, u_{t}(0, x) = u_{1}(x), \, \theta (0, x) = \theta_{0}(x), \, q(0, x) = q_{0}(x), \quad x \in \Omega \nonumber \\
		& u(t, x) = \theta(t, x) = 0, \quad (t, x) \in (0, \infty) \times \partial \Omega \nonumber
	\end{align}
	for $p > 2$, where a nonlinear source term is competing with the damping
	caused by the heat conduction and established a local existence result. He
	also showed that solutions with negative initial energy blow up in finite
	time. The blow-up result was then improved by Messaoudi and Said-Houari \cite{MeSa2004} 
	to accommodate certain solutions with positive initial energy.
	
	In the present work, we are concerned with 
	\begin{align}
		& \rho_{1} \varphi_{tt} - \sigma (\varphi_{x}, \psi)_{x} + \mu \varphi_{t} = 0, \quad (t, x) \in (0, \infty) \times (0, L),  \nonumber \\
		& \rho_{2} \psi_{tt} - b \psi_{xx} + k(\varphi_{x} + \psi) + \beta \theta_{x} = 0, \quad (t, x) \in (0, \infty) \times (0, L),  \nonumber \\
		& \rho_{3} \theta_{t} + \gamma q_{x} + \delta \psi_{tx} = 0, \quad (t, x) \in (0, \infty ) \times (0, L),  
		\label{TIMOSHENKO_CATTANEO_NONLINEAR_DAMPED} \\
		& \tau_{0} q_{t} + q + \kappa \theta_{x} = 0, \quad (t, x) \in (0, \infty) \times (0,L),  \nonumber
	\end{align}
	where $\varphi = \varphi(t, x)$ is the displacement vector, $\psi = \psi(t,x)$
	is the rotation angle of the filament, $\theta = \theta(t, x)$ is the
	temperature difference, $q = q(t, x)$ is the heat flux vector, $\rho_{1}$, 
	$\rho_{2}$, $\rho_{3}$, $b$, $k$, $\gamma$, $\delta$, $\kappa$, $\mu$, 
	$\tau_{0}$ are positive constants. The nonlinear function $\sigma$ is
	assumed to be sufficiently smooth and satisfy 
	\begin{equation}
		\sigma_{\varphi_{x}}(0, 0) = \sigma_{\psi}(0, 0) = k \nonumber
	\end{equation}
	and 
	\begin{equation}
		\sigma_{\varphi_{x} \varphi_{x}}(0, 0) = \sigma_{\varphi_{x}\psi}(0, 0) = \sigma_{\psi \psi} = 0.  \nonumber
	\end{equation}
	This system models the transverse vibration of a beam subject to the heat
	conduction given by Cattaneo's law instead of the usual Fourier's one. We
	should note here that dissipative effects of heat conduction induced by
	Cattaneo's law are usualy weaker than those induced by Fourier's law (an
	opposite effect was observed though in \cite{Ir2006}). This justifies the
	presence of the extra damping term in the first equation of (\ref{TIMOSHENKO_CATTANEO_NONLINEAR_DAMPED}). 
	In fact if $\mu = 0$, Fern\'{a}ndez Sare and Racke \cite{FeRa2007} have proved recently 
	that (\ref{TIMOSHENKO_CATTANEO_NONLINEAR_DAMPED}) is no longer exponentially stable
	even in the case of equal propagation speed ($\rho_{1}/\rho_{2} = k/b)$.
	Moreover, they showed that this ''unexpected`` phenomenon (the loss of
	exponential stability) takes place even in the presence of a viscoelastic
	damping in the second equation of (\ref{TIMOSHENKO_CATTANEO_NONLINEAR_DAMPED}). 
	If $\mu > 0$, but $\beta = 0$, one can also prove with the aid of semigroup theory (cf. \cite{MuRa2002}, Section 4) 
	that the system is not exponential stable independent of the relation between coefficients. 
	Our aim is to show that the presence of frictional damping $\mu \varphi_{t}$ in the first
	equation of (\ref{TIMOSHENKO_CATTANEO_NONLINEAR_DAMPED}) will drive the
	system to stability in an exponential rate independent of the wave speeds in linear and nonlinear cases.
	
	The structure of the paper is as follows. In section \ref{SECTION_LINEAR_1}, we discuss the well-posedness 
	and exponential stability of the linearized problem for $\varphi = \psi = q = 0$ on the boundary. 
	In section \ref{SECTION_LINEAR_2}, we establish the same result for $\varphi_{x} = \psi = q = 0$ on the boundary. 
	In section \ref{SECTION_NONLINEAR_1}, we study the nonlinear system subject to the boundary conditions $\varphi_{x} = \psi = q =0$, 
	show the global unique solvability and exponential stability for small initial data.

\section{Linear exponential stability --- $\varphi = \psi = q = 0$} \label{SECTION_LINEAR_1} 
	For the sake of technical convenience, by scaling the system (\ref{TIMOSHENKO_CATTANEO_NONLINEAR_DAMPED}), 
	we transform it to an equivalent form 
	\begin{align}
		&\rho_1 \varphi_{tt} - \sigma(\varphi_x, \psi)_x + \mu \varphi_t = 0, \quad (t, x) \in (0, \infty) \times (0, L),  \nonumber \\
		&\rho_2 \psi_{tt} - b \psi_{xx} + k(\varphi_x + \psi) + \gamma \theta_x = 0, \quad (t, x) \in (0, \infty) \times (0, L),  \nonumber \\
		&\rho_3 \theta_t + \kappa q_x + \gamma \psi_{tx} = 0, \quad (t, x) \in (0,
		\infty) \times (0, L),  \label{TIMOSHENKO_CATTANEO_NONLINEAR_DAMPED_EQUIV} \\
		&\tau_0 q_t + \delta q + \kappa \theta_x = 0, \quad (t, x) \in (0, \infty) \times (0, L),  \nonumber
	\end{align}
	with some other constants and the nonlinear function $\sigma$ still satisfying (possibly for a new $k$) 
	\begin{equation}
		\sigma_{\varphi_x}(0, 0) = \sigma_{\psi}(0, 0) = k \label{SIGMA_ASSUMPTION_1}
	\end{equation}
	and 
	\begin{equation}
		\sigma_{\varphi_x \varphi_x}(0, 0) = \sigma_{\varphi_x \psi}(0, 0) = \sigma_{\psi \psi} = 0.  \label{SIGMA_ASSUMPTION_2}
	\end{equation}
	
	In this section, we consider the linearization of (\ref {TIMOSHENKO_CATTANEO_NONLINEAR_DAMPED_EQUIV}) given by 
	\begin{align}
		&\rho_1 \varphi_{tt} - k(\varphi_x + \psi)_x + \mu \varphi_t = 0, \quad (t, x) \in (0, \infty) \times (0, L),  \nonumber \\
		&\rho_2 \psi_{tt} - b \psi_{xx} + k(\varphi_x + \psi) + \gamma \theta_x = 0, \quad (t, x) \in (0, \infty) \times (0, L),  \nonumber \\
		&\rho_3 \theta_t + \kappa q_x + \gamma \psi_{tx} = 0, \quad (t, x) \in (0, \infty) \times (0, L),  \label{TIMOSHENKO_CATTANEO_LINEAR_DAMPED_EQUIV} \\
		&\tau_0 q_t + \delta q + \kappa \theta_x = 0, \quad (t, x) \in (0, \infty) \times (0, L),  \nonumber
	\end{align}
	completed by the following boundary and initial conditions 
	\begin{align}
		\varphi(t, 0) &= \varphi(t, L) = \psi(t, 0) = \psi(t, L) = q(t, 0) = q(t, L) = 0 \text{ in } (0, \infty),
		\label{TIMOSHENKO_CATTANEO_LINEAR_DAMPED_EQUIV_BC} \\
		\varphi(0, \cdot) &= \varphi_0, \varphi_t(0, \cdot) = \varphi_1, \psi(0, \cdot) = \psi_0, \psi_t(0, \cdot) = \psi_1,  \nonumber \\
		\theta(0, \cdot) &= \theta_0, q(0, \cdot) = q_0. \label{TIMOSHENKO_CATTANEO_LINEAR_DAMPED_EQUIV_IC}
	\end{align}
	
	We present a brief discussion of the well-posedness, and the semigroup formulation of 
	(\ref{TIMOSHENKO_CATTANEO_LINEAR_DAMPED_EQUIV})--(\ref{TIMOSHENKO_CATTANEO_LINEAR_DAMPED_EQUIV_IC}). 
	For this purpose, we set $V:=(\varphi ,\varphi_{t},\psi ,\psi_{t},\theta ,q)^{t}$ and observe that $V$ satisfies 
	\begin{equation}
		\left\{ 
		\begin{array}{c}
			V_{t} = AV \\ 
			V(0) = V_{0}
		\end{array}
		\right., \label{TIMOSHENKO_CATTANEO_LINEAR_DAMPED_EQUIV_SEMIGROUP_FORMULATION}
	\end{equation}
	where $V_{0} := (\varphi_{0}, \varphi_{1}, \psi_{0}, \psi_{1},\theta_{0},q_{0})^{t}$ and $A$ is the differential operator 
	\begin{equation}
		A = \left( 
		\begin{array}{cccccc}
			0 & 1 & 0 & 0 & 0 & 0 \\ 
			\frac{k}{\rho_{1}}\partial_{x}^{2} & -\frac{\mu }{\rho_{1}} & \frac{k}{\rho_{1}}\partial_{x} & 0 & 0 & 0 \\ 
			0 & 0 & 0 & 1 & 0 & 0 \\ 
			-\frac{k}{\rho_{2}}\partial_{x} & 0 & \frac{b}{\rho_{2}} \partial_{x}^{2} -
			\frac{k}{\rho_{2}} & 0 & -\frac{\gamma}{\rho_{2}}\partial_{x} & 0 \\ 
			0 & 0 & 0 & -\frac{\gamma}{\rho_{3}}\partial_{x} & 0 & -\frac{\kappa}{\rho_{2}}\partial_{x} \\ 
			0 & 0 & 0 & 0 & -\frac{\kappa}{\tau_{0}}\partial_{x} & -\frac{\delta}{\tau_{0}}
		\end{array}
		\right).  \nonumber
	\end{equation}
	
	The energy space 
	\begin{equation}
		\mathcal{H} := H^1_0((0, L)) \times L^2((0, L)) \times H^1_0((0, L)) \times L^2((0, L)) \times L^2((0, L)) \times L^2((0, L)) \nonumber
	\end{equation}
	is a Hilbert space with respect to the inner product 
	\begin{align}
		\langle V, W \rangle_\mathcal{H} &= 
		\phantom{+} \rho_1 \langle V^1, W^1\rangle_{L^2((0, L))} + \rho_2 \langle V^4, W^4\rangle_{L^2((0, L))} \nonumber \\
		&\phantom{=} +b \langle V^3_x, W^3_x\rangle_{L^2((0, L))} + k \langle V^1_x + V^3, W^1_x + W^3\rangle_{L^2((0, L))}  \nonumber \\
		&\phantom{=} +\rho_3 \langle V^5, W^5\rangle_{L^2((0, L))} + \tau_0 \langle V^6, W^6\rangle  \nonumber
	\end{align}
	for all $V, W \in \mathcal{H}$. The domain of $A$ is then 
	\begin{align}
		D(A) = \{V \in \mathcal{H} \,|\, &V^1, V^3 \in H^2((0, L)) \cap H^1_0((0, L)), V^2, V^3 \in H^1_0((0, L))  \nonumber \\
		&V^5, V^6 \in H^1_0((0, L)), V^5_x \in H^1_0((0, L))\}.  \nonumber
	\end{align}
	
	It is easy to show according to \cite{Ra2002} the validness of	
	\begin{Lemma}
		The operator $A$ has the following properties:
		
		\begin{enumerate}
		\item  $\overline{D(A)}=\mathcal{H}$ and $A$ is closed;
		
		\item  $A$ is dissipative;
		
		\item  $D(A)=D(A^{\ast})$.
		\end{enumerate}
	\end{Lemma}
	
	Now, by the virtue of the Hille-Yosida theorem, we have the following result.
	
	\begin{Theorem}
		$A$ generates a $C_{0}$-semigroup of contractions $\{e^{At}\}_{t\geq 0}$. 
		If $V_{0} \in D(A)$, the unique solution $V \in C^{1}([0, \infty), \mathcal{H}) \cap C^{0}([0, \infty), D(A))$ 
		to (\ref{TIMOSHENKO_CATTANEO_LINEAR_DAMPED_EQUIV_SEMIGROUP_FORMULATION}) is given by 
		$V(t) = e^{At} V_{0}$. If $V_{0}\in D(A^{n})$ for $n \in \mathbb{N}$, then $V \in C^{0}([0, \infty), D(A^{n}))$.
	\end{Theorem}
	
	Our next aim is to obtain an exponential stability result for the energy
	functional $E(t) = E(t; \varphi, \psi, \theta, q)$ given by 
	\begin{align}
		E(t; \varphi, \psi, \theta, q) = \frac{1}{2} \int^L_0 (\rho_1 \varphi_t^2 +
		\rho_2 \psi_t^2 + b \psi_x^2 + k(\varphi_x + \psi)^2 + \rho_3 \theta^2 + \tau_0 q^2) \mathrm{d}x.  \nonumber
	\end{align}
	
	We formulate and prove the following theorem.
	
	\begin{Theorem}
		Let $(\varphi, \psi, \theta, q)$ be the unique solution to 
		(\ref{TIMOSHENKO_CATTANEO_LINEAR_DAMPED_EQUIV})--(\ref{TIMOSHENKO_CATTANEO_LINEAR_DAMPED_EQUIV_IC}). 
		Then, there exist two positive constants $C$ and $\alpha$, independent of $t$ and the initial data, such that 
		\begin{equation}
			E(t; \varphi, \psi, \bar{\theta}, q) \leq CE(0; \varphi, \psi, \bar{\theta}, q) e^{-2\alpha t} \text{ for all } t \geq 0, \nonumber
		\end{equation}
		where $\bar{\theta}(t, x) = \theta (t,x) - \frac{1}{L} \int_{0}^{L} \theta_0(s) \mathrm{d}s$.
	\end{Theorem}
	
	\begin{Proof}
		To show the exponential stability of the energy functional, we use the
		Lyapunov's method, i.e. we construct a Lyapunov functional $\mathcal{L}$ satisfying 
		\[
			\beta_{1}E(t) \leq \mathcal{L}(t) \leq \beta_{2}E(t), \quad t \geq 0 
		\]
		for positive constants $\beta_{1}$, $\beta_{2}$ and 
		\[
			\frac{\mathrm{d}}{\mathrm{d}t}\mathcal{L}(t)\leq -2\alpha \mathcal{L} (t),\quad t\geq 0 
		\]
		for some $\alpha > 0$. This will be achieved by a careful choice of multiplicators.
		
		Multiplying in $L^{2}((0,L))$ the first equation in (\ref{TIMOSHENKO_CATTANEO_LINEAR_DAMPED_EQUIV}) 
		by $\varphi_{t}$, the second by $\psi_{t}$, the third by $\theta$ and the fourth by $q$ and partially integrating, we obtain 
		\begin{equation}
			\frac{\mathrm{d}}{\mathrm{d}t}E(t) = -\mu \int_{0}^{L} \varphi_{t}^{2} \mathrm{d}x - \delta \int_{0}^{L}q^{2} \mathrm{d}x.
		\end{equation}
		
		As in \cite{MuRa2002}, let $w$ be a solution to 
		\[
			-w_{xx}=\psi_{x},\quad w(0)=w(L)=0 
		\]
		and let 
		\[
			I_{1} := \int_{0}^{L} \left(\rho_{2} \psi_{t} \psi + \rho_{1} \varphi_{t} w 
			- \frac{\gamma \tau_{0}}{\kappa} \psi q\right) \mathrm{d}x. 
		\]
		
		Then, we obtain taking into account the second equation in (\ref{TIMOSHENKO_CATTANEO_LINEAR_DAMPED_EQUIV}) 
		\begin{align}
		\frac{\mathrm{d}}{\mathrm{d}t}& \int_{0}^{L}\rho_{2}\psi_{t}\psi \mathrm{d}%
		x=\rho_{2}\int_{0}^{L}\left( \psi_{t}^{2}+\psi_{tt}\psi \right) \mathrm{d}%
		x  \nonumber \\
		& =\rho_{2}\int_{0}^{L}\psi_{t}^{2}\mathrm{d}x+b\int_{0}^{L}\psi_{xx}\psi 
		\mathrm{d}x-k\int_{0}^{L}(\varphi_{x}+\psi )\psi \mathrm{d}x-\gamma
		\int_{0}^{L}\theta_{x}\psi \mathrm{d}x.  \nonumber
		\end{align}
		
		Further, we get using the first and the fourth equations in (\ref{TIMOSHENKO_CATTANEO_LINEAR_DAMPED_EQUIV}) 
		\begin{align}
			\frac{\mathrm{d}}{\mathrm{d}t} &\int_{0}^{L} \rho_{1} \varphi_{t} w \mathrm{d}x =
			\rho_{1}\int_{0}^{L} \left(\varphi_{tt} w + \varphi_{t} w_{t}\right) \mathrm{d}x  \nonumber \\
			& = - k \int_{0}^{L} \varphi \psi_{x} \mathrm{d} x + k \int_{0}^{L} w_{x}^{2} \mathrm{d}x
			- \mu \int_{0}^{L} \varphi_{t} w \mathrm{d}x + \rho_{1} \int_{0}^{L} \varphi_{t}w_{t} \mathrm{d}x,  \nonumber \\
			\frac{\mathrm{d}}{\mathrm{d}t} & \int_{0}^{L} - \frac{\gamma \tau_{0}}{\kappa} \psi q\mathrm{d}x
			= - \frac{\gamma \tau_{0}}{\kappa } \int_{0}^{L} \psi_{t}q \mathrm{d}x
			+ \frac{\gamma}{\kappa} \int_{0}^{L} \psi (\delta q + \kappa \theta_{x})\mathrm{d}x  \nonumber \\
			& = -\frac{\gamma \tau_{0}}{\kappa} \int_{0}^{L} \psi_{t} q \mathrm{d}x 
			+ \frac{\gamma \delta}{\kappa} \int_{0}^{L}\psi q\mathrm{d}x + \gamma \int_{0}^{L} \theta_{x} \psi \mathrm{d}x.  \nonumber
		\end{align}
		
		By using the above inequalities, we find 
		\begin{align}
			\frac{\mathrm{d}}{\mathrm{d}t} I_{1}& = \rho_{2} \int_{0}^{L} \psi_{t}^{2} \mathrm{d}x - b \int_{0}^{L} \psi_{x}^{2} \mathrm{d}x 
			- k \int_{0}^{L} \psi^{2} \mathrm{d}x+k\int_{0}^{L}w_{x}^{2}\mathrm{d}x  \nonumber \\
			& \phantom{=} - \mu \int_{0}^{L} \varphi_{t} w \mathrm{d}x + \rho_{1} \int_{0}^{L} \varphi_{t} w_{t} \mathrm{d}x
			- \frac{\gamma \tau_{0}}{\kappa} \int_{0}^{L} \psi_{t} q \mathrm{d}x 
			+ \frac{\gamma \delta}{\kappa} \int_{0}^{L}\psi q\mathrm{d}x.  \nonumber
		\end{align}
		
		Observing 
		\begin{equation}
			\int_{0}^{L} w_{x}^{2} \mathrm{d}x \leq \int_{0}^{L} \psi^{2} \mathrm{d}x \leq c \int_{0}^{L} \psi_{x}^{2} \mathrm{d}x,
		\end{equation}
		with the Poincar\'{e} constant $c = \frac{L^{2}}{\pi^{2}} > 0$, we conclude using the Young's inequality 
		\begin{align}
			\frac{\mathrm{d}}{\mathrm{d}t} I_{1} 
			\leq & \phantom{+} \rho_{2} \int_{0}^{L}\psi_{t}^{2} \mathrm{d}x - b \int_{0}^{L} \psi_{x}^{2} \mathrm{d}x
			-k\int_{0}^{L}\psi^{2} \mathrm{d}x + k\int_{0}^{L} \psi^{2}\mathrm{d}x \nonumber \\
			& +\frac{\mu}{2} \int_{0}^{L} \left(\varepsilon_{1}w^{2} + \frac{1}{\varepsilon_{1}}\varphi_{t}^{2}\right) \mathrm{d}x
			+ \frac{\rho_{1}}{2} \int_{0}^{L} \left( \varepsilon_{1} w_{t}^{2} 
			+ \frac{1}{\varepsilon_{1}} \varphi_{t}^{2} \right) \mathrm{d}x  \nonumber \\
			& +\frac{\gamma \tau_{0}}{2\kappa} \int_{0}^{L} \left(\varepsilon_{1} \psi_{t}^{2}+\frac{1}{\varepsilon_{1}} q^{2}\right)\mathrm{d}x
			+\frac{\gamma \delta}{2\kappa} \int_{0}^{L} \left(\varepsilon_{1} \psi^{2} + \frac{1}{\varepsilon_{1}}q^{2}\right) \mathrm{d}x  \nonumber \\
			\leq &-\left[b - \frac{\varepsilon_{1}}{2} \left(\mu c^{2} + \frac{\delta \gamma c}{\kappa }\right) \right] 
			\int_{0}^{L}\psi_{x}^{2} \mathrm{d} x + \left[\rho_{2} + \frac{\varepsilon_{1}}{2} 
			\left(\rho_{1}c + \frac{\gamma \tau_{0}}{\kappa }\right) \right] \int_{0}^{L}\psi_{t}^{2}\mathrm{d}x  \nonumber \\
			& +\frac{1}{2\varepsilon_{1}}\left( \mu +\rho_{1}\right) \int_{0}^{L}\varphi_{t}^{2}\mathrm{d}x
			+ \frac{1}{2\varepsilon_{1}}\left(\frac{\gamma \tau_{0}}{\kappa } + \frac{\delta \gamma}{\kappa}\right)
			\int_{0}^{L}q^{2}\mathrm{d}x.  \label{I1_ESTIMATE}
		\end{align}
		for some $\varepsilon_{1} > 0$.
		
		Next, we consider the functional $I_{2}$ given by 
		\[
			I_{2}:=\rho_{1}\int_{0}^{L}\varphi_{t}\varphi \mathrm{d}x. 
		\]
		
		It easily follows that 
		\begin{align}
			\frac{\mathrm{d}}{\mathrm{d}t} I_{2} &= \rho_{1} \int_{0}^{L} \varphi_{tt} \varphi \mathrm{d}x
			+\rho_{1} \int_{0}^{L} \varphi_{t}^{2} \mathrm{d}x \nonumber \\
			& =\int_{0}^{L} k(\varphi_{x} + \psi)_{x} \varphi \mathrm{d}x 
			- \mu \int_{0}^{L}\varphi_{t}\varphi \mathrm{d}x + \rho_{1}\int_{0}^{L}\varphi_{t}^{2}\mathrm{d}x  \nonumber \\
			& =-k\int_{0}^{L}\varphi_{x}^{2} \mathrm{d}x + k\int_{0}^{L} \psi_{x} \varphi \mathrm{d}x
			- \mu \int_{0}^{L}\varphi_{t}\varphi \mathrm{d}x + \rho_{1}\int_{0}^{L}\varphi_{t}^{2}\mathrm{d}x,  \nonumber
		\end{align}
		which can be estimated by 
		\begin{align}
			\frac{\mathrm{d}}{\mathrm{d}t}I_{2}& \leq -k \int_{0}^{L} \varphi_{x}^{2} \mathrm{d}x
			+\frac{k}{2} \int_{0}^{L} \left(\varepsilon_{2} \varphi^{2} + \frac{1}{\varepsilon_{2}}\psi_{x}^{2}\right) \mathrm{d}x  \nonumber \\
			& \phantom{=} + \frac{\mu}{2} \int_{0}^{L} \left(\varepsilon_{2}\varphi^{2} 
			+ \frac{1}{\varepsilon_{2}}\varphi_{t}^{2} \right) \mathrm{d}x
			+\rho_{1}\int_{0}^{L}\varphi_{t}^{2}\mathrm{d}x  \nonumber
		\end{align}
		\begin{align}
			& \leq - \left(k - \frac{\varepsilon_{2}c}{2}(k + \mu)\right)
			\int_{0}^{L} \varphi_{x}^{2} \mathrm{d}x + \frac{k}{2\varepsilon_{2}} \int_{0}^{L}\psi_{x}^{2}\mathrm{d}x  \nonumber \\
			& \phantom{=} +\left(\frac{\mu}{2\varepsilon_{2}} + \rho_{1}\right) \int_{0}^{L} \varphi_{t}^{2} \mathrm{d}x  \label{I2_ESTIMATE}
		\end{align}
		for some $\varepsilon_{2}>0$.
		
		Next we consider a functional $I_{3}$ defined by 
		\begin{equation}
			I_{3} := N_{1} I_{1} + I_{2}  \nonumber
		\end{equation}
		for some $N_{1} > 0$ and, combining (\ref{I1_ESTIMATE}) and (\ref{I2_ESTIMATE}), arrive at 
		\begin{align}
			\frac{\mathrm{d}}{\mathrm{d}t}I_{3} \leq& -\left[N_{1} \left(b - \frac{\varepsilon_{1}}{2} 
			\left(\mu c^{2} +\frac{\delta \gamma c}{\kappa}\right)\right) 
			- \frac{k}{2\varepsilon_{2}}\right] \int_{0}^{L} \psi_{x}^{2} \mathrm{d}x  \nonumber \\
			&-\left(k - \frac{\varepsilon_{2}c}{2}(k + \mu)\right) 
			\int_{0}^{L} \varphi_{x}^{2} \mathrm{d}x + N_{1} \left[\rho_{2} + \frac{\varepsilon_{1}}{2}
			\left(\rho_{1}c+\frac{\gamma \tau_{0}}{\kappa }\right) \right] \int_{0}^{L}\psi_{t}^{2}\mathrm{d}x  \nonumber \\
			& +\left[ N_{1}\frac{1}{2 \varepsilon_{1}}\left(\mu + \rho_{1}\right) + \left(\frac{\mu}{2\varepsilon_{2}} + \rho_{1}\right) \right]
			\int_{0}^{L}\varphi_{t}^{2} \mathrm{d}x  \nonumber \\
			& +N_{1}\frac{1}{2\varepsilon_{1}} \left(\frac{\gamma \tau_{0}}{\kappa} 
			+ \frac{\delta \gamma}{\kappa}\right) \int_{0}^{L}q^{2} \mathrm{d}x. \label{I3_ESTIMATE}
		\end{align}
		
		At this point, we introduce 
		\begin{equation}
			\bar{\theta}(t, x) = \theta(t, x) - \frac{1}{L} \int_{0}^{L} \theta_{0}(x) \mathrm{d}x.  \nonumber
		\end{equation}
		
		One can easily verify that $(\varphi, \psi, \bar{\theta}, q)$ satisfies system (\ref{TIMOSHENKO_CATTANEO_LINEAR_DAMPED_EQUIV}). 
		Moreover, one can apply the Poincar\'{e} inequality to $\bar{\theta}$ 
		\begin{equation}
			\int_{0}^{L}\bar{\theta}^{2}(t, x) \mathrm{d}x \leq c \int_{0}^{L}\bar{\theta}_{x}^{2}(t,x)\mathrm{d}x,  \nonumber
		\end{equation}
		since $\int_{0}^{L} \bar{\theta}(t, x) \mathrm{d}x = 0$ for all $t \geq 0$. 
		Until the end of this chapter, we shall work with $\bar{\theta}$ but denote it with $\theta$.
		
		In order to obtain a negative term of $\int_{0}^{L} \psi_{t}^{2} \mathrm{d}x$, we introduce, as in \cite{MuRa2002}, 
		the following functional 
		\[
			I_{4}(t) := \rho_{2} \rho_{3} \int_{0}^{L} \left(\int_{0}^{x} \theta(t, y) \mathrm{d}y\right) \psi_{t}(t, x)\mathrm{d}x, 
		\]
		and find 
		\begin{align}
			\frac{\mathrm{d}}{\mathrm{d}t} I_{4} 
			=& \int_{0}^{L} \left(\int_{0}^{x} \rho_{3}\theta_{t}\mathrm{d}y\right) \rho_{2}\psi_{t} \mathrm{d}x
			+\int_{0}^{L} \left(\int_{0}^{x} \rho_{3} \theta \mathrm{d}y\right) \rho_{2} \psi_{tt}\mathrm{d}x  \nonumber \\
			=& -\int_{0}^{L} \left(\int_{0}^{x}\kappa q_{x} + \gamma \psi_{tx} \mathrm{d}y\right) \rho_{2}\psi_{t}\mathrm{d}x  \nonumber \\
			& +\int_{0}^{L} \left(\int_{0}^{x} \rho_{3} \theta \mathrm{d}y\right) 
			(b\psi_{xx} - k(\varphi_{x} + \psi) - \gamma \theta_{x})\mathrm{d}x  \nonumber
		\end{align}
		\begin{align}
			=& -\gamma \rho_{2} \int_{0}^{L} \psi_{t}^{2} \mathrm{d}x - \rho_{2} \kappa \int_{0}^{L} q\psi_{t} \mathrm{d}x
			- b\rho_{3}\int_{0}^{L} \theta \psi_{x} \mathrm{d}x  \nonumber \\
			& + k\rho_{3} \int_{0}^{L} \theta \varphi \mathrm{d}x 
			- k\rho_{3}\int_{0}^{L} \left(\int_{0}^{x}\theta \mathrm{d}y\right) \psi \mathrm{d} x
			+ \gamma \rho_{3} \int_{0}^{L} \theta ^{2} \mathrm{d}x.  \nonumber
		\end{align}
		
		This can be estimated as follows 
		\begin{align}
			\frac{\mathrm{d}}{\mathrm{d}t} I_{4} \leq & -\gamma \rho_{2} \int_{0}^{L} \psi_{t}^{2} \mathrm{d}x
			+ \frac{\rho_{2}\kappa}{2} \int_{0}^{L} \left(\varepsilon_{4} \psi_{t}^{2} + \frac{1}{\varepsilon_{4}}q^{2}\right) \mathrm{d}x
			+\frac{b\rho_{3}}{2} \int_{0}^{L} \varepsilon_{4}^{\prime} \psi_{x}^{2} \nonumber \\
			& +\frac{1}{\varepsilon_{4}^{\prime}} \theta ^{2}\mathrm{d}x 
			+ \frac{k\rho_{3}}{2} \int_{0}^{L} \left(\varepsilon_{4}^{\prime} \varphi^{2}
			+\frac{1}{\varepsilon_{4}^{\prime}} \theta^{2} \right) \mathrm{d}x
			+\frac{k\rho_{3}}{2} \int_{0}^{L} \varepsilon_{4}^{\prime} \psi^{2} \mathrm{d}x  \nonumber \\
			& +\frac{1}{\varepsilon_{4}^{\prime}} \left(\int_{0}^{x} \theta \mathrm{d}y\right)^{2} \mathrm{d}x
			+ \gamma \rho_{3} \int_{0}^{L} \theta^{2} \mathrm{d}x  \nonumber \\
			=& \left[-\gamma \rho_{2} + \frac{\varepsilon_{4} \rho_{2}\kappa}{2}\right] 
			\int_{0}^{L} \psi_{t}^{2} \mathrm{d}x + \left(\frac{\varepsilon_{4}^{\prime} \rho_{3}}{2}(b + kc)\right) \int_{0}^{L}\psi_{x}^{2}\mathrm{d}x \nonumber \\
			& +\frac{\varepsilon_{4}^{\prime} k\rho_{3}c}{2} \int_{0}^{L}\varphi_{x}^{2}\mathrm{d}x
			+ \left(\gamma \rho_{3} + \frac{\rho_{3}}{2\varepsilon_{4}^{\prime}}(b+k+kc)\right) \int_{0}^{L}\theta ^{2}\mathrm{d}x  \nonumber\\
			& +\frac{\rho_{2}\kappa}{2\varepsilon_{4}} \int_{0}^{L} q^{2} \mathrm{d}x \label{I4_ESTIMATE}
		\end{align}
		for arbitrary positive $\varepsilon_{4}$ and $\varepsilon_{4}^{\prime}$.
		
		Finally, we set 
		\[
			I_{5}(t) := -\tau_{0} \rho_{3} \int_{0}^{L} q(t, x) \left(\int_{0}^{x} \theta(t, y) \mathrm{d}y \right) \mathrm{d}x
		\]
		and observe 
		\begin{align}
			\frac{\mathrm{d}}{\mathrm{d}t}I_{5}(t) &= 
			-\rho_{3} \int_{0}^{L} \tau_{0} q_{t}\left(\int_{0}^{x} \theta \mathrm{d}y\right) \mathrm{d}x 
			-\tau_{0} \int_{0}^{L} q\left(\int_{0}^{x} \rho_{3} \theta_{t}\mathrm{d}y\right) \mathrm{d}x  \nonumber \\
			&= -\rho_{3} \int_{0}^{L}(-\delta q - \kappa \theta_{x})
			\left(\int_{0}^{x}\theta \mathrm{d}y\right) \mathrm{d}x \nonumber \\
			&\phantom{=} -\tau_{0}\int_{0}^{L} q\left(\int_{0}^{x} - \kappa q_{x} - \gamma \psi_{tx}\mathrm{d}y\right) \mathrm{d}x  \nonumber \\
			&= \rho_{3} \delta \int_{0}^{L} q \left(\int_{0}^{x} \theta \mathrm{d}y\right) 
			\mathrm{d}x + \rho_{3} \kappa \int_{0}^{L} \theta_{x} \left(\int_{0}^{x}\theta \mathrm{d}y\right) \mathrm{d}x  \nonumber \\
			&\phantom{=} +\tau_{0} \kappa \int_{0}^{L} q \left(\int_{0}^{x} q_{x} \mathrm{d}y\right) \mathrm{d}x
			+\tau_{0} \gamma \int_{0}^{L} q \left(\int_{0}^{x} \psi_{tx} \mathrm{d}y\right) \mathrm{d}x  \nonumber \\
			&= \frac{\rho_{3}\delta}{2} \int_{0}^{L} \left(\varepsilon_{5}\left(\int_{0}^{x} \theta^{2} \mathrm{d}y\right)^{2}
			+ \frac{1}{\varepsilon_{5}} q^{2}\right) \mathrm{d}x - \rho_{3}\kappa \int_{0}^{L} \theta^{2} \mathrm{d}x  \nonumber \\
			&\phantom{=}+ \tau_{0}\kappa \int_{0}^{L} q^{2} \mathrm{d}x 
			+ \frac{\tau_{0}\gamma}{2} \int_{0}^{L}\varepsilon_{5}^{\prime} \psi_{t}^{2}
			+\frac{1}{\varepsilon_{5}^{\prime}} q^{2}\mathrm{d}x  \nonumber
		\end{align}
		\begin{align}
			& \leq \left(-\rho_{3} \kappa + \frac{\varepsilon_{5}\rho_{3}\delta c}{2}\right) \int_{0}^{L} \theta^{2}\mathrm{d}x
			+ \frac{\varepsilon_{5}^{\prime} \tau_{0} \gamma}{2} \int_{0}^{L} \psi_{t}^{2} \mathrm{d}x  \nonumber \\
			&\phantom{=}+ \left(\tau_{0}\kappa + \frac{\rho_{3}\delta}{2\varepsilon_{5}} 
			+ \frac{\tau_{0}\gamma}{2\varepsilon_{5}^{\prime}}\right) \int_{0}^{L}q^{2} \mathrm{d}x  \label{I5_ESTIMATE}
		\end{align}
		for positive $\varepsilon_{5}$ and $\varepsilon_{5}^{\prime }$
		
		For $N,N_{4},N_{5}>0$, we can define an auxiliary functional $\mathcal{F}(t)$ by 
		\[
			\mathcal{F}(t) := NE + I_{3} + N_{4} I_{4} + N_{5} I_{5}.
		\]
		
		From (\ref{I3_ESTIMATE}), (\ref{I4_ESTIMATE}) and (\ref{I5_ESTIMATE}), we have then 
		\begin{align}
			\frac{\mathrm{d}}{\mathrm{d}t} \mathcal{F}(t) \leq& -C_{\psi_{x}} \int_{0}^{L} \psi_{x}^{2} \mathrm{d}x
			- C_{\varphi_{x}} \int_{0}^{L} \varphi_{x}^{2} \mathrm{d}x - C_{\psi_{t}} \int_{0}^{L} \psi _{t}^{2}\mathrm{d}x  \nonumber \\
			&- C_{\theta} \int_{0}^{L} \theta^{2} \mathrm{d}x - C_{\varphi_{t}} \int_{0}^{L} \varphi_{t}^{2} \mathrm{d}x 
			-C_{q}\int_{0}^{L} q^{2} \mathrm{d}x,  \label{LYAPUNOV_ESTIMATE_1}
		\end{align}
		where 
		\begin{align}
			C_{\psi_{x}} &= \left[N_{1} \left(b - \frac{\varepsilon_{1}}{2} 
			\left(\mu c^{2} + \frac{\delta \gamma c}{\kappa}\right) \right) 
			- \frac{k}{2\varepsilon_{2}} - N_{4} \frac{\varepsilon_{4}^{\prime}}{2} \rho_{3}(b + kc)\right],  \nonumber \\
			C_{\varphi_{x}} &= \left[\left(k - \frac{\varepsilon_{2}}{2}c(k + \mu)\right) 
			-N_{4} \frac{\varepsilon_{4}^{\prime}}{2} k\rho_{3} c\right],  \nonumber \\
			C_{\psi_{t}} &= \left[N_{4} \left(\gamma \rho_{2} -\frac{\varepsilon_{4} \rho_{2}\kappa}{2}\right) 
			-N_{1} \left(\rho_{2} + \frac{\varepsilon_{1}}{2} \left(\rho_{1} c + \frac{\gamma \tau_{0}}{\kappa}\right)\right) 
			-N_{5} \frac{\varepsilon_{5}^{\prime} \tau_{0} \gamma}{2}\right],  \nonumber \\
			C_{\theta} &= \left[N_{5} \left(\rho_{3} \kappa - \frac{\varepsilon_{5} 
			\rho_{3} \delta c}{2}\right) - N_{4} \left(\gamma \rho_{3} 
			+ \frac{\rho_{3}}{2\varepsilon_{4}^{\prime}} (b + k + kc)\right)\right],  \nonumber \\
			C_{\varphi_{t}} &= \left[N \mu - N_{1} \frac{1}{2 \varepsilon_{1}}
			(\mu + \rho_{1}) - \left(\frac{\mu}{2\varepsilon_{2}} + \rho_{1}\right) \right], \nonumber \\
			C_{q} &= \left[N - N_{1} \frac{1}{2\varepsilon_{1}} 
			\left(\frac{\gamma \tau_{0}}{\kappa} + \frac{\delta \gamma}{\kappa}\right) - N_{4} \frac{\rho_{2} \kappa}{2\varepsilon_{4}} 
			- N_{5} \left(\tau_{0} \kappa + \frac{\rho_{3}\delta}{2\varepsilon_{5}} 
			+ \frac{\tau_{0}\gamma}{2\varepsilon_{5}^{\prime}}\right)\right] .  \nonumber
		\end{align}
		
		Choosing $\varepsilon_{1}$, $\varepsilon_{2}$, $\varepsilon_{4}$, $\varepsilon_{5}$ sufficiently small, 
		then $N_{1}$ and $N_{4}$ sufficiently large, $\varepsilon_{4}^{\prime}$ sufficiently small, 
		$N_{5}$ sufficiently large, $\varepsilon_{5}^{\prime}$ sufficiently small and finally $N$ sufficiently large, we can assure that 
		\begin{align*}
			\varepsilon_{1} &< \frac{2b\kappa}{\mu \kappa c^{2} + \delta \gamma c}, \quad
			\varepsilon_{2} < \frac{2k}{c(k + \mu)}, \quad 
			\varepsilon_{4} < \frac{2\gamma}{\kappa}, \quad 
			\varepsilon_{5} < \frac{2\kappa}{\delta c}, \\
			N_{1} &> \frac{k}{2\varepsilon_{2} \left(b - \frac{\varepsilon_{1}}{2} \left(\mu c^{2} 
			+ \frac{\delta \gamma c}{\kappa }\right) \right)}, \\
			N_{4} &> \frac{N_{1} \left(\rho_{2} + \frac{\varepsilon_{1}}{2} \left(\rho_{1} c
			+ \frac{\gamma \tau_{0}}{\kappa}\right) \right)}{\gamma \rho_{2}- \frac{\varepsilon_{4}\rho_{2}\kappa }{2}}, \\
			\varepsilon_{4}^{\prime}& < \min\left\{\frac{2N_{1} \left(b - \frac{\varepsilon_{1}}{2} \left(\mu c^{2}
			+\frac{\delta \gamma c}{\kappa }\right)\right)}{N_{4}\rho_{3}(b + kc)}, 
			\frac{2\left(k - \frac{\varepsilon_{2}}{2} c(k + \mu)\right)}{N_{4} k \rho_{3}c}\right\},
		\end{align*}
		\begin{align*}
			N_{5} &> \frac{N_{4} \left(\gamma \rho_{3} + \frac{\rho_{3}}{2\varepsilon_{4}^{\prime }} 
			\left(b + k + kc\right)\right)}{\rho_{3}\kappa -\frac{\varepsilon_{5} \rho_{3} \delta c}{2}}, \\
			\varepsilon_{5}^{\prime} &< \frac{2\left(N_{4} \left(\gamma \rho_{2} - \frac{\varepsilon_{4}\rho_{2}\kappa}{2}\right) 
			- N_{1}\left( \rho_{2} + \frac{\varepsilon_{1}}{2} \left(\rho_{1} c 
			+ \frac{\gamma \tau_{0}}{\kappa} \right)\right)\right)}{N_{5} \tau_{0}\gamma} \\
			N &> \max\Bigg\{\frac{N_{1} \frac{1}{2 \varepsilon_{1}}(\mu + \rho_{1}) + \left(\frac{\mu}{2\varepsilon_{2}}
			+ \rho_{1}\right)}{\mu}, \\
			&\phantom{> \max\{}N_{1} \frac{1}{2\varepsilon_{1}} \left(\frac{\gamma \tau_{0}}{\kappa} + \frac{\delta \gamma}{\kappa}\right) 
			+ N_{4}\frac{\rho_{2}\kappa}{2\varepsilon_{4}} + N_{5} \left(\tau_{0}\kappa + \frac{\rho_{3}\delta}{2\varepsilon_{5}}
			+\frac{\tau_{0}\gamma}{2\varepsilon_{5}^{\prime}}\right) \Bigg \}.
		\end{align*}
		
		Having fixed the constants as above, we find that all the terms on the right-hand side of (\ref{LYAPUNOV_ESTIMATE_1}) are negative.
		
		Now, we have to estimate $\frac{\mathrm{d}}{\mathrm{d}t}\mathcal{F}(t)$ versus $-d_{2}E(t)$ for a $d_{2} > 0$. 
		By letting $C := \frac{1}{2}\min\{C_{\psi_{x}},C_{\varphi_{x}}\}$, we conclude from (\ref{LYAPUNOV_ESTIMATE_1}) that 
		\begin{align}
			\frac{\mathrm{d}}{\mathrm{d}t} \mathcal{F}(t) &\leq 
			\underbrace{-C\int_{0}^{L}\psi_{x}^{2}\mathrm{d}x}_{\leq - \frac{C}{c}\int_{0}^{L}\psi^{2}\mathrm{d}x}
			- C\int_{0}^{L}\varphi_{x}^{2}\mathrm{d}x -(C_{\psi_{x}}-C) \int_{0}^{L}\psi_{x}^{2}\mathrm{d}x  \nonumber \\
			&-C_{\psi_{t}} \int_{0}^{L} \psi_{t}^{2} \mathrm{d}x 
			- C_{\theta} \int_{0}^{L}\theta ^{2}\mathrm{d}x-C_{\varphi_{t}}\int_{0}^{L}\varphi _{t}^{2}\mathrm{d}x
			-C_{q}\int_{0}^{L}q^{2}\mathrm{d}x  \nonumber \\
			&\leq - \min\left\{C, \frac{C}{c}\right\} \int_{0}^{L} \underbrace{\left(\varphi_{x}^{2} 
			+ \psi ^{2}\right)}_{\geq \frac{1}{2}(\varphi_{x}+\psi)^{2}}\mathrm{d}x 
			- (C_{\psi_{x}}-C)\int_{0}^{L}\psi_{x}^{2}\mathrm{d}x \nonumber \\
			& -C_{\psi_{t}} \int_{0}^{L} \psi_{t}^{2} \mathrm{d}x - C_{\theta} \int_{0}^{L} \theta^{2} \mathrm{d}x
			- C_{\varphi_{t}} \int_{0}^{L}\varphi_{t}^{2} \mathrm{d}x - C_{q} \int_{0}^{L} q^{2} \mathrm{d}x  \nonumber \\
			& \leq -C_{\varphi_{t}}\int_{0}^{L}\varphi_{t}^{2}\mathrm{d}x-C_{\psi_{t}}
			\int_{0}^{L}\psi_{t}^{2}\mathrm{d}x - (C_{\psi_{x}} - C) \int_{0}^{L} \psi_{x}^{2}\mathrm{d}x  \nonumber \\
			& -\frac{\min \left\{C, \frac{C}{c}\right\}}{2} \int_{0}^{L}(\varphi_{x}
			+\psi )^{2}\mathrm{d}x - C_{\theta} \int_{0}^{L} \theta^{2} \mathrm{d} x - C_{q} \int_{0}^{L} q^{2}\mathrm{d}x  \nonumber \\
			& \leq -d_{1} \int_{0}^{L}(\varphi_{t}^{2} + \psi_{t}^{2} + \psi_{x}^{2} + (\varphi_{x}+\psi )^{2} 
			+\theta^{2} + q^{2}) \mathrm{d}x. \label{LYAPUNOV_ESTIMATE_2}
		\end{align}
		with 
		\begin{equation}
			d_{1} : =\min \left\{ C_{\varphi_{t}},C_{\psi_{t}},(C_{\psi_{x}}-C),
			\frac{\min \left\{C, \frac{C}{c}\right\}}{2},C_{\theta}, C_{q}\right\} .
		\end{equation}
		
		For $d_{2} := \frac{2d_{1}}{\max\{\rho_{1}, \rho_{2}, b, k, \rho_{3}, \tau_{0}\}}$, we can therefore estimate 
		\[
			\frac{\mathrm{d}}{\mathrm{d}t}\mathcal{F}(t) \leq -d_{2}E(t).
		\]
		
		Finally, we consider the functional $H(t) := I_{3} + N_{4} I_{4} + N_{5} I_{5}$ and show for this 
		\[
			|H(t)| \leq CE(t), \quad C > 0.
		\]
		
		By using the trivial relation 
		\[
			\int_{0}^{L} \varphi^{2} \mathrm{d}x \leq 2c \int_{0}^{L}(\varphi_{x} 
			+ \psi)^{2}\mathrm{d}x+2c^{2}\int_{0}^{L} \psi_{x}^{2} \mathrm{d}x
		\]
		with the Poincar\'{e} constant $c=\frac{L^{2}}{\pi ^{2}}$ we arive at 
		\begin{align}
			|H(t)| &= \left|N_{1} I_{1} + I_{2} + N_{4} I_{4} + N_{5} I_{5}\right| 
			\leq N_{1} |I_{1}| + |I_{2}| + N_{4} |N_{4}| + N_{5} |I_{5}|  \nonumber \\
			&= N_{1} \left|\int_{0}^{L} \left(\rho_{2} \psi_{t} \psi 
			+ \rho_{1}\varphi_{t}w - \frac{\gamma \tau_{0}}{\kappa}\psi q\right) \mathrm{d}x\right| 
			+ \rho_{1}\left| \int_{0}^{L} \varphi_{t} \varphi \mathrm{d}x\right|   \nonumber \\
			&+ N_{4} \rho_{2} \rho_{3} \left|\int_{0}^{L} \left(\int_{0}^{x}\theta (t,x) \mathrm{d}y\right) \psi_{t}(t,x)\mathrm{d}x\right| 
			+ N_{5}\tau_{0}\rho_{3} \left|\int_{0}^{L} q\left(\int_{0}^{x} \theta \mathrm{d}y\right) \mathrm{d}x\right|  \nonumber \\
			&\leq N_{1} \Bigg(\frac{\rho_{2}}{2} \int_{0}^{L} \psi_{t}^{2} \mathrm{d}x+
			\frac{\rho_{2}c}{2}\int_{0}^{L} \psi_{x}^{2} \mathrm{d}x 
			+\frac{\rho_{1}}{2} \int_{0}^{L}\varphi_{t}^{2} \mathrm{d}x 
			+\frac{\rho_{1}c^{2}}{2} \int_{0}^{L}\psi_{x}^{2} \mathrm{d}x  \nonumber \\
			&+ \frac{\gamma \tau_{0}c}{2\kappa} \int_{0}^{L} \psi_{x}^{2} \mathrm{d}x
			+ \frac{\gamma \tau_{0}}{2\kappa} \int_{0}^{L} q^{2} \mathrm{d}x\Bigg) 
			+\frac{\rho_{1}}{2} \left(\int_{0}^{L}\varphi_{t}^{2} \mathrm{d} x + \int_{0}^{L}\varphi^{2} \mathrm{d}x\right)  \nonumber \\
			&+ \frac{\rho_{2} \rho_{3} N_{4}}{2} \left(c \int_{0}^{L}\theta^{2} \mathrm{d}x
			+ \int_{0}^{L} \psi_{t}^{2} \mathrm{d}x \right) + \frac{\tau_{0} \rho_{3} N_{5}}{2}
			\left( \int_{0}^{L}q^{2}\mathrm{d}x+c\int_{0}^{L}\theta ^{2}\mathrm{d} x\right)   \nonumber \\
			&\leq \hat{C}_{\varphi_{t}} \int_{0}^{L} \varphi_{t}^{2} + \hat{C}_{\psi_{t}} \int_{0}^{L}\psi_{t}^{2} \mathrm{d}x
			+ \hat{C}_{\varphi_{x}} \int_{0}^{L}\psi_{x}^{2} \mathrm{d}x  \nonumber \\
			&+ \hat{C}_{\varphi_{x} + \psi} \int_{0}^{L}(\varphi_{x} + \psi)^{2} \mathrm{d}x
			+ \hat{C}_{\theta} \int_{0}^{L}\theta^{2} \mathrm{d}x + \hat{C}_{q} \int_{0}^{L} q^{2}\mathrm{d}x,  \label{H_ESTIMATE}
		\end{align}
		where the constants are determined as follows 
		\begin{align}
			\hat{C}_{\varphi_{t}} &:= \frac{1}{2}\left(N_{1} \rho_{1} + \rho_{1}\right), \quad 
			\hat{C}_{\psi_{t}} := \frac{1}{2} \left(N_{1} \rho_{2} + \rho_{2} \rho_{3}N_{4}\right) ,  \nonumber \\
			\hat{C}_{\psi_{x}} &:= \frac{1}{2} \left(N_{1} \rho_{2} c + N_{1} \rho_{1} c^{2} 
			+ \frac{N_{1}\tau_{0}c}{\kappa }+2\rho_{1}c^{2}\right) ,  \nonumber \\
			\hat{C}_{\varphi_{x} + \psi} &:= \rho_{1} c, \, 
			\hat{C}_{\theta} := \frac{1}{2} \left(N_{4} \rho_{2} \rho_{3} c +N_{5} \rho_{3} \tau_{0} c\right), \,
			\hat{C}_{q} := \frac{1}{2} \left(\frac{N_{1} \gamma \tau_{0}}{\kappa}
			+ N_{5} \rho_{3} \tau_{0}\right).  \nonumber
		\end{align}
		
		According to (\ref{H_ESTIMATE}) we have $|H(t)|\leq \hat{C}E(t)$ for 
		\[
			\hat{C} := \frac{\max\left\{\hat{C}_{\varphi_{t}}, \hat{C}_{\psi_{t}},
			\hat{C}_{\psi_{x}}, \hat{C}_{\varphi_{x} + \psi}, \hat{C}_{\theta}, \hat{C}_{q}\right\}}
			{\min\left\{\rho_{1}, \rho_{2}, b, k, \rho_{3}, \tau_{0}\right\}}.
		\]
		
		Taking finally $\hat{N} > \max\{N, \hat{C}\}$ and defining a Lyapunov functional 
		\begin{equation}
			\mathcal{L}(t) := \hat{N} E + H(t) = \hat{N} E + I_{3} + N_{4} I_{4} + N_{5} I_{5}, \label{LYAPUNOV_FUNCTIONAL_FINAL}
		\end{equation}
		we obtain, on the one hand, 
		\begin{equation}
			\beta_{1} E(t)\leq \mathcal{L}(t) \leq \beta_{2}E(t) \label{LYAPUNOV_FUNCTIONAL_EQUIVALENCE}
		\end{equation}
		for $\beta_{1} := \hat{N} - \hat{C} > 0$, $\beta_{2} := \hat{N} + \hat{C} > 0$, 
		on the other hand, we know that 
		\[
			\frac{\mathrm{d}}{\mathrm{d}t} \mathcal{L}(t) \leq -d_{2}E(t) \leq - \frac{d_{2}}{\beta_{2}} \mathcal{L}(t).
		\]
		
		By using the Gronwall's lemma, we conclude for $\alpha := \frac{d_{2}}{2\beta_{2}}$ that 
		\[
			\mathcal{L}(t)\leq e^{-2\alpha t}0).
		\]
		
		Eventually, (\ref{LYAPUNOV_FUNCTIONAL_EQUIVALENCE}) yields 
		\[
			E(t)\leq Ce^{-2\alpha t}E(0)
		\]
		with $C:=\frac{\beta_{2}}{\beta_{1}}$.
	\end{Proof}

\section{Linear exponential stability --- $\varphi_x = \psi = q = 0$} \label{SECTION_LINEAR_2} 
	The second set of boundary conditions we are going to study in this paper is 
	\begin{equation}
		\varphi_x(t, 0) = \varphi_x(t, L) = \psi(t, 0) = \psi(t, L) = q(t, 0) = q(t, L) = 0 \text{ in } (0, \infty).
		\label{TIMOSHENKO_CATTANEO_LINEAR_DAMPED_EQUIV_BC_1} \\
	\end{equation}
	
	Here, we consider the initial boundary value problem (\ref{TIMOSHENKO_CATTANEO_LINEAR_DAMPED_EQUIV}), 
	(\ref{TIMOSHENKO_CATTANEO_LINEAR_DAMPED_EQUIV_IC}), 
	(\ref{TIMOSHENKO_CATTANEO_LINEAR_DAMPED_EQUIV_BC_1}). 
	We will present a semigroup formulation of this problem, show the exponential stability of the
	associated semigroup and make estimates on higher energies. This will enable
	us to prove global existence and exponential stability also in nonlinear
	settings.
	
	Let 
	\begin{align}
		L^2_\ast((0, L)) &= \big\{u \in L^2((0, L)) \,\big|\, \int^L_0 u(x) \mathrm{d}x = 0\big\},  \nonumber \\
		H^1_\ast((0, L)) &= \big\{u \in H^1((0, L)) \,\big|\, \int^L_0 u(x) \mathrm{d}x = 0\big\}.  \nonumber
	\end{align}
	
	We introduce a Hilbert space 
	\begin{equation}
		\mathcal{H} := H_{\ast}^{1}((0,L)) \times L_{\ast}^{2}((0,L))\times
		H_{0}^{1}((0,L)) \times L^{2}((0,L)) \times L_{\ast}^{2}((0,L)) \times L^{2}((0,L)) \nonumber
	\end{equation}
	equipped with the inner product 
	\begin{align}
		\langle V,W\rangle_{\mathcal{H}} &= \phantom{+} \rho_{1} \langle V^{1}, W^{1}\rangle_{L^{2}((0,L))}
		+ \rho_{2} \langle V^{4}, W^{4} \rangle_{L^{2}((0,L))}  \nonumber \\
		&\phantom{=} + b\langle V_{x}^{3}, W_{x}^{3} \rangle_{L^{2}((0, L))} 
		+ k\langle V_{x}^{1} + V^{3}, W_{x}^{1} + W^{3}\rangle_{L^{2}((0, L))}  \nonumber \\
		&\phantom{=} + \rho_{3} \langle V^{5}, W^{5}\rangle_{L^{2}((0, L))} + \tau_{0}\langle V^{6},W^{6}\rangle_{L^{2}((0,L))}.  \nonumber
	\end{align}
	
	Let the operator $A$ be formally defined as in section \ref{SECTION_LINEAR_1}
	with the domain 
	\begin{align}
		D(A) = \{V \in \mathcal{H} \,|\, &V^1 \in H^2((0, L)), V^1_x \in H^1_0((0, L)), V^2 \in H^1_\ast((0, L)),  \nonumber \\
		&V^3 \in H^2((0, L)), V^4 \in H^1_0((0, L)),  \nonumber \\
		&V^5 \in H^1_\ast((0, L)), V^6 \in H^1_0((0, L))\}.  \nonumber
	\end{align}
	
	Setting $V := (\varphi, \varphi_{t}, \psi, \psi_{t}, \theta ,q)^{t}$, we
	observe that $V$ satisfies 
	\begin{equation}
		\left\{ 
		\begin{array}{c}
			V_{t}=AV \\ 
			V(0)=V_{0}
		\end{array}
		\right. ,
		\label{TIMOSHENKO_CATTANEO_LINEAR_DAMPED_EQUIV_SEMIGROUP_FORMULATION_1}
	\end{equation}
	where $V_{0} := (\varphi_{0}, \varphi_{1}, \psi_{0}, \psi_{1}, \theta_{0}, q_{0})^{t}$.
	
	By assuring that $A$ satisfies the conditions of the Hille-Yosida theorem, we can easily get
	
	\begin{Theorem}
		$A$ generates a $C_{0}$-semigroup of contractions $\{e^{At}\}_{t \geq 0}$. 
		If $V_{0}\in D(A)$, the the unique solution $V\in C^{1}([0, \infty),\mathcal{H}) \cap C^{0}([0, \infty), D(A))$ 
		to (\ref{TIMOSHENKO_CATTANEO_LINEAR_DAMPED_EQUIV_SEMIGROUP_FORMULATION_1}) is given by $V(t) = e^{At}V_{0}$. 
		If $V_{0}\in D(A^{n})$ for $n \in \mathbb{N}$, then $V \in C^{0}([0, \infty), D(A^{n}))$.
	\end{Theorem}
	
	Moreover, we can show that the Lyapunov functional (\ref{LYAPUNOV_FUNCTIONAL_FINAL}) constructed in
	section \ref{SECTION_LINEAR_1} is also a Lyapunov functional for (\ref{TIMOSHENKO_CATTANEO_LINEAR_DAMPED_EQUIV_SEMIGROUP_FORMULATION_1}).
	Observing for the energy $E(t)$ of the unique solution $(\varphi, \psi, \theta, q)$ that 
	\begin{equation}
		E(t) = \frac{1}{2} \|V \|_{\mathcal{H}}^{2}  \nonumber
	\end{equation}
	holds independent of $t$, we obtain the exponential stability of the associated semigroup $\{e^{At}\}_{t\geq 0}$.
	
	\begin{Theorem} \label{THEOREM_EXPONENTIAL_STABILITY_LINEAR_2} 
		The semigroup $\{e^{At}\}_{t \geq 0}$ associated with $A$ is exponential stable, i.e. 
		\begin{equation}
			\exists c_{1} > 0 \quad \forall t \geq 0 \quad \forall V_{0} \in \mathcal{H} : 
			\quad \|e^{At} V_{0} \|_{\mathcal{H}} 
			\leq c_{1} e^{-\alpha t} \|V_{0}\|_{\mathcal{H}}.  \label{EXPONENTIAL_STABILITY_SEMIGROUP}
		\end{equation}
	\end{Theorem}
	
	Similar to \cite{MuRa2002}, we observe that if $V_{0} \in D(A)$, we can estimate $AV(t)$ in the same way as $V(t)$ is estimated 
	in (\ref{EXPONENTIAL_STABILITY_SEMIGROUP}), implying in its turn using the structure of $A$ 
	that $(V_{x}^{1},V_{x}^{2},V_{x}^{3},V_{x}^{4},V_{x}^{5},V_{x}^{6})$ can be estimated in the norm of $\mathcal{H}$, 
	hence, one can estimate $((\varphi_{x})_{x}, (\varphi_{t})_{x}, (\psi_{x})_{x}, (\psi_{t})_{x}, \theta_{x}, q_{x})^{t}$ in $L^{2}((0,L))^{6}$.
	
	We define for $s \in \mathbb{N}$ the Hilbert space 
	\begin{equation}
		\mathcal{H}_s := (H^{s} \times H^{s-1} \times H^{s} \times H^{s-1} \times H^{s-1} \times H^{s-1})((0, L))  \nonumber
	\end{equation}
	with natural norm Sobolev norm for its component. Using the consideration above, we can therefore estimate 
	\begin{equation}
		\|V(t)\|_{\mathcal{H}_s} \leq c_s \|V_0\|_{\mathcal{H}_s} e^{-\alpha t}. \label{ESTIMATE_HIGHER_ENERGIES}
	\end{equation}
	$c_s$ denotes here a positive constant, being independent of $V_0$ and $t$.

\section{Nonlinear exponential stability} \label{SECTION_NONLINEAR_1}
	In this section, we study the nonlinear system 
	\begin{align}
		&\rho_1 \varphi_{tt} - \sigma(\varphi_x, \psi)_x + \mu \varphi_t = 0, \quad (t, x) \in (0, \infty) \times (0, L),  \nonumber \\
		&\rho_2 \psi_{tt} - b \psi_{xx} + k(\varphi_x + \psi) + \gamma \theta_x = 0, \quad (t, x) \in (0, \infty) \times (0, L),  \nonumber \\
		&\rho_3 \theta_t + \kappa q_x + \gamma \psi_{tx} = 0, \quad (t, x) \in (0, \infty) \times (0, L),  \label{TIMOSHENKO_CATTANEO_NONLINEAR_DAMPED_EQUIV_1} \\
		&\tau_0 q_t + \delta q + \kappa \theta_x = 0, \quad (t, x) \in (0, \infty) \times (0, L),  \nonumber
	\end{align}
	completed by the boundary 
	\begin{align}
		\varphi(t, 0) &= \varphi(t, L) = \psi(t, 0) = \psi(t, L) = q(t, 0) = q(t, L) = 0 \text{ in } (0, \infty),
		\label{TIMOSHENKO_CATTANEO_NONLINEAR_DAMPED_EQUIV_BC_1}
	\end{align}
	and the initial conditions 
	\begin{align}
		\varphi(0, \cdot) &= \varphi_0, \quad \varphi_t(0, \cdot) = \varphi_1, \quad
		\psi(0, \cdot) = \psi_0, \quad \psi_t(0, \cdot) = \psi_1,  \nonumber \\
		\theta(0, \cdot) &= \theta_0, \quad q(0, \cdot) = q_0.
		\label{TIMOSHENKO_CATTANEO_NONLINEAR_DAMPED_EQUIV_IC_1}
	\end{align}
	
	As before, the constants $\rho_1$, $\rho_2$, $\rho_3$, $b$, $k$, $\gamma$, $\delta$, $\kappa$, $\mu$, $\tau_0$ are assumed to be positive. 
	The nonlinear function $\sigma$ is assumed to be sufficiently smooth and to satisfy 
	\begin{equation}
		\sigma_{\varphi_x}(0, 0) = \sigma_{\psi}(0, 0) = k \label{SIGMA_ASSUMPTION_1_1}
	\end{equation}
	and 
	\begin{equation}
		\sigma_{\varphi_x \varphi_x}(0, 0) = \sigma_{\varphi_x \psi}(0, 0) = \sigma_{\psi \psi} = 0.  \label{SIGMA_ASSUMPTION_2_1}
	\end{equation}
	
	To obtain a local well-posedness result, we have first to consider a
	corresponding non-homogeneous linear system 
	\begin{align}
		&\rho_1 \varphi_{tt} - \hat{\sigma}(t, x) \varphi_{xx} - \check{\sigma}(t, x) \psi_x + \mu \varphi_t = 0 \quad \text{ in } (0, \infty) \times (0, L),  \nonumber \\
		&\rho_2 \psi_{tt} - b\psi_{xx} + k(\varphi_x + \psi) + \gamma \theta_x = 0 \quad \text{ in } (0, \infty) \times (0, L),  \nonumber \\
		&\rho_3 \theta_t + \kappa q_x + \gamma \psi_{tx} = 0 \quad \text{ in } (0, \infty) \times (0, L),  \label{TIMOSHENKO_CATTANEO_LINEAT_NONHOMEGENEOUS} \\
		&\tau_0 q_t + \delta q + \kappa \theta_{x} = 0 \quad \text{ in } (0, \infty) \times (0, L)  \nonumber
	\end{align}
	together with the boundary conditions (\ref{TIMOSHENKO_CATTANEO_NONLINEAR_DAMPED_EQUIV_BC_1}) 
	and initial conditions (\ref{TIMOSHENKO_CATTANEO_NONLINEAR_DAMPED_EQUIV_IC_1}).
	
	The solvability of this system is established in the following theorem.
	
	\begin{Theorem} \label{LOCAL_LINEAR_EXISTENCE_THEOREM} 
		We assume for some $T>0$ that 
		\begin{align}
			&\hat{\sigma},\check{\sigma}\in C^{1}([0,T]\times [0,L]),  \nonumber \\
			&\hat{\sigma}_{tt},\hat{\sigma}_{tx},\hat{\sigma}_{xx},\check{\sigma}_{tt},
			\check{\sigma}_{tx},\check{\sigma}_{xx}\in L^{\infty}([0,T], L^{2}((0,L))). \nonumber
		\end{align}
		
		Let $\hat{\sigma}\geq s > 0$. The initial data may satisfy 
		\begin{align}
			\varphi_{0,x} &\in H^{2}((0, L)) \cap H_{0}^{1}((0, L)), \quad \varphi_{1, x} \in H_{0}^{1}((0, L)),  \nonumber \\
			\psi_{0} &\in H^{3}((0, L))\cap H_{0}^{1}((0, L)), \quad \psi_{1}\in H^{2}((0,L))\cap H_{0}^{1}((0,L)),  \nonumber \\
			\theta_{0} &\in H^{2}((0,L)),\quad q_{0}\in H^{2}((0,L))\cap H_{0}^{1}((0,L)).  \nonumber
		\end{align}
		
		Under the above conditions, the initial boundary problem (\ref{TIMOSHENKO_CATTANEO_LINEAT_NONHOMEGENEOUS}), 
		(\ref{TIMOSHENKO_CATTANEO_NONLINEAR_DAMPED_EQUIV_BC_1}), (\ref{TIMOSHENKO_CATTANEO_NONLINEAR_DAMPED_EQUIV_IC_1}), posesses a unique
		classical solution $(\varphi, \psi, \theta, q)$ such that 
		\begin{align}
			\varphi, \psi &\in C^{2}([0, T] \times [0, L]), \quad \theta, q \in C^{1}([0,T] \times [0, L]),  \nonumber \\
			\partial^{\alpha} \varphi, \partial^{\alpha} \psi &\in L^{\infty}([0,T], L^{2}((0, L))), \quad 1 \leq |\alpha| \leq 3,  \nonumber \\
			\partial^{\alpha} \theta, \partial^{\alpha}q &\in L^{\infty}([0, T], L^{2}((0, L))), \quad 0 \leq |\alpha| \leq 2  \nonumber
		\end{align}
		with $\partial^{\alpha} = \partial_{t}^{\alpha_{1}}\partial_{x}^{\alpha_{2}}$ 
		for $\alpha = (\alpha_{1}, \alpha_{2}) \in \mathbb{N}_{0}^{2}$.
	\end{Theorem}
	
	\begin{Proof}
		We present here a similar proof to that one of Slemrod in \cite{Sl1981}.
		Using Faedo-Galerkin method, we construct a sequence that converges to a
		solution of (\ref{TIMOSHENKO_CATTANEO_LINEAT_NONHOMEGENEOUS}), 
		(\ref{TIMOSHENKO_CATTANEO_NONLINEAR_DAMPED_EQUIV_BC_1}), (\ref{TIMOSHENKO_CATTANEO_NONLINEAR_DAMPED_EQUIV_IC_1}). 
		By using then a special a priori estimate, one obtains corresponding regularity of the solution.
		
		Letting $\lambda_{i} := i\pi/L$, $c_{i}(x) := \sqrt{2/L} \cos\lambda_{i}x$, 
		$s_{i}(x) := \sqrt{2/L} \sin \lambda_{i}x$, $i \in \mathbb{N}$, 
		we define $(\varphi_{m}(t), \psi_{m}(t), \theta_{m}(t), q_{m}(t))$ by 
		\begin{align}
			\varphi_{m}(t) &:= \sum_{i = 0}^{m} \Phi_{im}(t) c_{i}(x), \quad 
			\psi_{m}(t) := \sum_{i = 0}^{m} \Psi_{im}(t) s_{i}(x),  \nonumber \\
			\theta_{m}(t) &:= \sum_{i = 0}^{m} \Theta_{im}(t) c_{i}(x),\quad
			q_{m}(t) := \sum_{i = 0}^{m} Q_{im}(t) s_{i}(x),  \nonumber
		\end{align}
		where 
		\begin{align}
			\Phi_{im}(0) &= \int_{0}^{L} \varphi_{0}(x) c_{i}(x) \mathrm{d}x, \quad 
			\dot{\Phi}_{im}(0) = \int_{0}^{L}\varphi_{1}(x) s_{i}(x) \mathrm{d}x,  \nonumber \\
			\Psi_{im}(0) &= \int_{0}^{L} \psi_{0}(x) s_{i}(x) \mathrm{d}x,
			\quad \dot{\Psi}_{im}(0) = \int_{0}^{L} \psi_{1}(x) c_{i}(x) \mathrm{d}x,  \nonumber \\
			\Theta_{im}(0) &= \int_{0}^{L} \theta_{0}(x) c_{i}(x) \mathrm{d}x, \quad
			Q_{im}(0) = \int_{0}^{L} q_{0}(x) s_{i}(x) \mathrm{d}x.  \nonumber
		\end{align}
		
		Multiplying the equations in (\ref{TIMOSHENKO_CATTANEO_LINEAT_NONHOMEGENEOUS}) 
		in $L^{2}((0,L))$ by $c_{i}$, $s_{i}$, $c_{i}$ and $s_{i}$, respectively,
		we observe that the functions $\Phi_{im}$, $\Psi_{im}$, $\Theta_{im}$, $Q_{im}$ 
		satisfy a system of ordinary differential equations 
		\begin{align}
			\rho_{1} \ddot{\Phi}_{jm}(t) =& -\sum_{i=0}^{m} \Phi_{im}(t) \lambda_{i}^{2}\langle \hat{\sigma}(t,x) c_{i}(x), c_{j}(x)\rangle   \nonumber \\
			\phantom{=}& + \sum_{i=0}^{m} \lambda_{i} \Psi_{im} \langle \check{\sigma}(t, x)c_{i} (x), c_{j}(x)\rangle 
			-\mu \dot{\Phi}_{jm}(t)  \nonumber \\
			\rho_{3} \ddot{\Psi}_{jm}(t) =& -b\Psi_{jm}(t) \lambda_{j}^{2} + k(\Phi_{jm}(t)\lambda_{j} 
			- \Psi_{jm}(t)) + \gamma \Theta_{jm}(t) \lambda_{j}, \label{GALERKIN_EQUATIONS} \\
			\rho_{3} \dot{\Theta}_{jm}(t) =& -\kappa Q_{jm}(t) \lambda_{j} - \gamma \dot{\Psi}_{jm}(t)\lambda_{j},  \nonumber \\
			\tau_{0} \dot{Q}_{jm}(t) =& -Q_{jm}(t) + \kappa \Theta_{jm}(t)\lambda_{j} \nonumber
		\end{align}
		for $0\leq j\leq m$ and 
		\[
			\langle f,g\rangle =\langle f,g\rangle_{L^{2}((0,L))}=\int_{0}^{L}f(x)g(x)\mathrm{d}x.
		\]
		
		This system is always solvable and possesses a unique solution 
		\begin{equation}
			(\Phi_{jm}, \Psi_{jm}, \Theta_{jm}, Q_{jm})  \nonumber
		\end{equation}
		with $\Phi_{jm}, \Psi_{jm} \in C^{2}([0, T])$ and $\Theta_{jm}, Q_{jm} \in C^{1}([0, T])$.
		
		We define a total energy $\mathcal{E}$ by 
		\begin{align}
			\mathcal{E}(t) =& \phantom{+} E(t; \varphi, \psi, \theta, q)
			+ E(t; \varphi_{t}, \psi_{t}, \theta_{t}, q_{t}) + E(t; \varphi_{tt}, \psi_{tt}, \theta_{tt}, q_{tt})  \nonumber \\
			&+ E(t; \varphi_{x}, \psi_{x}, \theta_{x}, q_{x}) + E(t; \varphi_{tx}, \psi_{tx}, \theta_{tx}, q_{tx}),  \nonumber
		\end{align}
		where 
		\begin{equation}
			E(t; \phi_{m}, \psi_{m}, \theta_{m}, q_{m}) = \frac{1}{2} \int_{0}^{L}(\rho_{1} \varphi_{t}^{2}
			+ \hat{\sigma} \varphi_{x}^{2} + \rho_{2}\psi_{t}^{2} + b\psi_{x}^{2} + \rho_{3}\theta^{3} + \tau_{0}q^{2})(t, x)\mathrm{d}x.
			\nonumber
		\end{equation}
		
		By multiplying in $L^{2}((0, L))$ the equations in (\ref{GALERKIN_EQUATIONS})
		by $\dot{\Phi}_{jm}$, $\dot{\Psi}_{jm}$, $\Theta_{jm}$, $Q_{jm}$, then
		differentiating them once and twice with respect to $t$, multiplying them
		with $\ddot{\Phi}_{jm}$, $\ddot{\Psi}_{jm}$, $\dot{\Theta}_{jm}$, $\dot{Q}_{jm}$ 
		and $\dddot{\Phi}_{jm}$, $\dddot{\Psi}_{jm}$, $\ddot{\Theta}_{jm}$, $\ddot{Q}_{jm}$, 
		respectively, and summing up over $j = 1, \dots, m$, we obtain an energy equality of the form 
		\begin{equation}
			\frac{\text{\textrm{d}}}{\text{\textrm{d}}t} \mathcal{E}(t; \varphi_{m},\psi_{m}, \theta_{m}, q_{m})
			= F_{m}(\partial^{\alpha_{1}}\varphi, \partial^{\alpha_{2} \psi}, \partial^{\beta_{1}} \theta,
			\partial^{\beta_{2}}q)  \nonumber
		\end{equation}
		for $0\leq |\alpha_{1,2}|\leq 3,$ $0\leq |\beta_{1,2}|\leq 2.$
		
		Following the approach of Slemrod and obtaining higher order $x$ derivatives
		from differential equations, we can integrate the above equality with
		respect to $t$ and estimate 
		\begin{equation}
			\int_{0}^{t} F_{m}(\tau) \mathrm{d} \tau 
			\leq C \int_{0}^{t} \mathcal{E}(\tau; \varphi_{m}, \psi_{m}, \theta_{m}, q_{m}) \mathrm{d}\tau .  \nonumber
		\end{equation}
		
		Gronwall's inequality yields then $\mathcal{E}(t) \leq C\mathcal{E}(0) e^{Ct} \leq C$ for a generic constant $C > 0$.
		
		It follows that the sequence $\{(\varphi_{m}, \psi_{m}, \theta_{m}, q_{m})\}_m$
		has a convergent subsequence. By the virtue of usual Sobolev embedding
		theorems, we get necessary regularity of the solution.
		
		The solution is unique since our a priori estimate can be shown also for $(\varphi, \psi, \theta, q)$ 
		assuring the continuous dependence of the solution on the initial data. 
		By usual continuation arguments, the solution can be smoothly continued to a maximal open interval $[0, T)$.
	\end{Proof}
	
	Having proved the local linear existence theorem, we can obtain a local existence also in the nonlinear situation.
	
	\begin{Theorem} \label{LOCAL_NONLINEAR_EXISTENCE_THEOREM} 
		Consider the initial boundary value problem 
		(\ref{TIMOSHENKO_CATTANEO_NONLINEAR_DAMPED_EQUIV_1})---(\ref{TIMOSHENKO_CATTANEO_NONLINEAR_DAMPED_EQUIV_IC_1}). 
		Let $\sigma = \sigma(r, s)\in C^{3}(\mathbb{R} \times \mathbb{R})$ satisfy 
		\begin{align}
			&0 < r_{0}\leq \sigma_{r} \leq r_{1} < \infty \quad (r_{0}, r_{1} > 0), \label{SIGMA_NL_COND_1} \\
			&0 \leq |\sigma_{s}| \leq s_{0} < \infty \quad (s_{0} > 0). \label{SIGMA_NL_COND_1}
		\end{align}
		
		Let the initial data comply with 
		\begin{align}
			\varphi_{0, x} &\in H^{2}((0, L)) \cap H_{0}^{1}((0, L)), \quad \varphi_{1, x}\in H_{0}^{1}((0, L)),  \nonumber \\
			\psi_{0} &\in H^{3}((0, L)) \cap H_{0}^{1}((0, L)), \quad 
			\psi_{1} \in H^{2}((0, L)) \cap H_{0}^{1}((0, L)),  \nonumber \\
			\theta_{0} &\in H^{2}((0, L)), \quad q_{0} \in H^{2}((0, L)) \cap H_{0}^{1}((0, L)).  \nonumber
		\end{align}
		
		The problem (\ref{TIMOSHENKO_CATTANEO_NONLINEAR_DAMPED_EQUIV_1})---(\ref{TIMOSHENKO_CATTANEO_NONLINEAR_DAMPED_EQUIV_IC_1}) 
		has then a unique classical solution $(\varphi, \psi, \theta, q)$ with 
		\begin{align}
			\varphi,\psi & \in C^{2}([0, T)\times [0, L]),  \nonumber \\
			\theta, q &\in C^{1}([0, T)\times [0, L]),  \nonumber
		\end{align}
		defined on a maximal existence interval $[0, T)$, $T \leq \infty$ such that for all $t_{0} \in [0, T)$ 
		\begin{align}
			\partial^{\alpha} \varphi, \partial^{\alpha} \psi &\in L^{\infty}([0, t_{0}],L^{2}((0, L))), 
			\quad 1 \leq |\alpha| \leq 3,  \nonumber \\
			\partial^{\alpha} \theta, \partial^{\alpha} q &\in L^{\infty}([0, t_{0}], L^{2}((0, L))),
			\quad 0 \leq |\alpha| \leq 2  \nonumber
		\end{align}
		holds.
	\end{Theorem}
	
	\begin{Proof}
		The proof of the local existence is by now standard. For positive $M$, $T$, we
		define the space $X(M, T)$ to be a set of all functions $(\varphi, \psi, \theta, q)$ such that they satisfy 
		\begin{align}
			\varphi(0, \cdot) &= \varphi_{0}, \quad \psi(0, \cdot) = \psi_{0}, \quad
			\theta(0, \cdot) = \theta_{0}, \quad q(0, \cdot) = q_{0},  \nonumber \\
			\varphi_{t}(0, \cdot ) &= \varphi_{1}, \quad \psi_{t}(0, \cdot) = \psi_{1} \quad \text{in }(0, L),  \label{ABNLA_X} \\
			\varphi_{x}(t, 0) &= \varphi_{x}(t, L) = \psi(t, 0) = \psi(t, L) = q(t, 0) = q(t, L) = 0 
			\quad \text{ in }(0, \infty )  \label{RBNL_X}
		\end{align}
		and their generalized derivatives fulfil 
		\begin{align}
			\partial^{\alpha} \varphi, \partial^{\alpha} \psi &\in L^{\infty}([0, T], L^{2}((0, L))), 
			\quad 1 \leq |\alpha | \leq 3,  \nonumber \\
			\partial^{\alpha} \theta, \partial^{\alpha} q &\in L^{\infty}([0, T], L^{2}((0, L))),
			\quad 0 \leq |\alpha| \leq 2  \nonumber
		\end{align}
		and 
		\[
			\sup_{0 \leq t \leq T} \int_{0}^{L} \left(\sum_{|\alpha|=1}^{3}
			\left[(\partial^{\alpha} \varphi)^{2} + (\partial^{\alpha} \psi)^{2}\right]
			+ \sum_{|\alpha| = 0}^{2}\left[(\partial^{\alpha}\theta)^{2} 
			+(\partial^{\alpha} q)^{2}\right]\right) \mathrm{d}x\leq M^{2}.
		\]
		
		Let $(\bar{\varphi},\bar{\psi},\bar{\theta},\bar{q})\in X(M,T)$. Consider
		the linear initial boundary value problem 
		\begin{align}
			&\rho_{1} \varphi_{tt} - \sigma_{r}(\bar{\varphi}_{x}, \bar{\psi})\varphi_{xx} 
			- \sigma_{s}(\bar{\varphi}_{x}, \bar{\psi})\psi_{x} + \mu \varphi_{t} = 0,\nonumber \\
			&\rho_{2} \psi_{tt} - b\psi_{xx} + k(\varphi_{x} + \psi) + \gamma \theta_{x} = 0, \nonumber \\
			&\rho_{3} \theta_{t} + \kappa q_{x} + \gamma \psi_{tx} = 0, \label{TIMOSHENKO_S_MAPPING_SYSTEM} \\
			&\tau_{0} q_{t} + \delta q + \kappa \theta_{x} = 0  \nonumber
		\end{align}
		together with initial conditions (\ref{TIMOSHENKO_CATTANEO_NONLINEAR_DAMPED_EQUIV_IC_1}) 
		and boundary conditions (\ref{TIMOSHENKO_CATTANEO_NONLINEAR_DAMPED_EQUIV_BC_1}).
		
		We set 
		\begin{align}
			\hat{\sigma}(t,x) &= \sigma_{r}(\bar{\varphi}_{x}, \bar{\psi}),  \nonumber \\
			\check{\sigma}(t,x) &= \sigma_{s}(\bar{\varphi}_{x}, \bar{\psi}),  \nonumber
		\end{align}
		and observe that $\hat{\sigma}$, $\check{\sigma}$ and the initial data
		satisfy the assumptions of the local existence and uniqueness Theorem \ref{LOCAL_LINEAR_EXISTENCE_THEOREM}. 
		Therefore, this linear problem possesses a unique solution.
		
		We define an operator $S$ mapping $(\bar{\varphi}, \bar{\psi}, \bar{\theta}, \bar{q})\in X(M,T)$ 
		to the solution of (\ref{TIMOSHENKO_S_MAPPING_SYSTEM}),
		i.e. $S(\bar{\varphi},\bar{\psi},\bar{\theta}, \bar{q}) = (\varphi, \psi, \theta, q)$.
		
		With standard techniques, we can show that $S$ maps the space $X(M, T)$ into
		itself if $M$ is sufficiently big and $T$ sufficiently small. 
		Following the approach of Slemrod, we show that $S$ is a contraction for sufficiently small $T$. 
		As $X(M,T)$ is a closed subset of the metric space 
		\begin{equation}
			Y = \{(\varphi ,\psi ,\theta ,q) \,|\, 
			\varphi_{t}, \varphi_{x}, \psi_{t}, \psi_{t}, \theta, q \in L^{\infty }([0, T],L^{2}((0, L)))\}  \nonumber
		\end{equation}
		equipped with a distance function 
		\begin{align}
			\rho\left((\varphi, \psi, \theta, q), (\bar{\varphi}, \bar{\psi}, \bar{\theta}, \bar{q})\right) 
			:= \sup_{0 \leq t \leq T} &\int_{0}^{L} \Big[(\varphi_{t} - \bar{\varphi}_{t})^{2}
			+(\varphi_{x} - \bar{\varphi}_{x})^{2} + (\psi_{t} - \bar{\psi}_{t})^{2}  \nonumber \\
			+& (\psi_{x} - \bar{\psi}_{x})^{2} + (\theta - \bar{\theta})^{2} + (q - \bar{q})^{2}\Big] \mathrm{d}x,  \nonumber
		\end{align}
		the Banach mapping theorem is applicable to $S$ and yields a unique solution in $X(M, T)$ having the asserted regularity.
	\end{Proof}
	
	To be able to handle the nonlinear problem globally, we need a local
	existence theorem with higher regularity. This one can be proved in the same
	way as Theorem \ref{LOCAL_NONLINEAR_EXISTENCE_THEOREM}.
	
	\begin{Theorem}
		\label{LOCAL_NONLINEAR_EXISTENCE_THEOREM_HIGHER_REGULARITY} Consider the initial boundary value problem
		(\ref{TIMOSHENKO_CATTANEO_NONLINEAR_DAMPED_EQUIV_1})---(\ref{TIMOSHENKO_CATTANEO_NONLINEAR_DAMPED_EQUIV_IC_1}). 
		Let $\sigma = \sigma(r, s) \in C^{4}(\mathbb{R} \times \mathbb{R})$ satisfy 
		\begin{align}
			&0 < r_{0} \leq \sigma_{r} \leq r_{1} < \infty \quad (r_{0}, r_{1} > 0),  \nonumber \\
			&0 \leq |\sigma_{s} |\leq s_{0} < \infty \quad (s_{0} > 0).  \nonumber
		\end{align}
		
		Let the assumptions of Theorem \ref{LOCAL_NONLINEAR_EXISTENCE_THEOREM} be satisfied. Moreover, let us assume 
		\[
			\varphi_{0,xxxx}, \psi_{0,xxxx}, \varphi_{1,xxx}, \psi_{1, xxx}, \theta_{0, xxx}, q_{0, xxx} \in L^{2}((0, L))
		\]
		and 
		\begin{align}
			\partial_{t}^{2} \varphi(0, \cdot), \partial_{t}^{2} \psi(0, \cdot ) &\in H^{2}((0, L)), \quad 
			\partial_{t}^{2} \varphi_{x}(0, \cdot), \partial_{t}^{2} \psi(0, \cdot) \in H_{0}^{1}((0, L))  \nonumber \\
			\partial_{t} \theta(0, \cdot), \partial_{t} q(0, \cdot) &\in H^{2}((0, L)), \quad 
			\partial_{t} q(0, \cdot) \in H_{0}^{1}((0, L)).  \nonumber
		\end{align}
		
		Then, (\ref{TIMOSHENKO_CATTANEO_NONLINEAR_DAMPED_EQUIV_1})---(\ref{TIMOSHENKO_CATTANEO_NONLINEAR_DAMPED_EQUIV_IC_1}) 
		possesses a unique classical solution $(\varphi ,\psi ,\theta ,q)$ satisfying 
		\begin{align}
			\varphi,\psi &\in C^{3}([0, T) \times [0, L]),  \nonumber \\
			\theta, q &\in C^{2}([0, T) \times [0, L]),  \nonumber
		\end{align}
		being defined in a maximal existence interval $[0, T)$, $T \leq \infty$ such that for all $t_{0} \in [0,T)$ 
		\begin{align}
			\partial^{\alpha} \varphi, \partial^{\alpha} \psi &\in L^{\infty}([0, t_{0}],L^{2}((0, L))), 
			\quad 1\leq |\alpha| \leq 4,  \nonumber \\
			\partial^{\alpha} \theta, \partial^{\alpha} q &\in L^{\infty}([0, t_{0}], L^{2}((0,L))), 
			\quad 0 \leq |\alpha| \leq 3  \nonumber
		\end{align}
		holds. Moreover, this interval coincides with that one from Theorem \ref{LOCAL_NONLINEAR_EXISTENCE_THEOREM}.
	\end{Theorem}
	
	\begin{Remark}
		Our conjecture is that in analogy to thermoelastic equations one can prove a
		more general existence theorem by getting bigger regularity of the solution
		under the same regularity assumptions as in Theorem \ref{LOCAL_NONLINEAR_EXISTENCE_THEOREM} for initial data (cf. \cite{JiRa2000}).
		
		This technique dates back to Kato and is based on a general notion of a CD-system coming from the semigroup theory.
	\end{Remark}
	
	For the proof of global solvability and exponential stability, we rewrite
	the problem (\ref{TIMOSHENKO_CATTANEO_NONLINEAR_DAMPED_EQUIV_1})---(\ref{TIMOSHENKO_CATTANEO_NONLINEAR_DAMPED_EQUIV_IC_1}) 
	as a nonlinear evolution problem.
	
	Letting $V = (\varphi, \varphi_{t}, \psi, \psi_{t}, \theta, q)^{t}$ and defining a linear differential operator 
	$A : D(A) \subset \mathcal {H} \to \mathcal{H}$ in the same manner as in section \ref{SECTION_LINEAR_2}, we obtain 
	\begin{equation}
		\left\{ 
		\begin{array}{c}
			V_{t} = AV + F(V, V_{x}) \\ 
			V(0) = V_{0}
		\end{array}
		\right. 
		\label{TIMOSHENKO_CATTANEO_NONLINEAR_DAMPED_EQUIV_SEMIGROUP_FORMULATION}
	\end{equation}
	with a nonlinear mapping $F$ being defined by 
	\begin{align}
		F(V,V_{x}) &= (0, \sigma_{\varphi_{x}}(\varphi_{x}, \psi) \varphi_{xx}
		- k \varphi_{xx} + \sigma_{\psi}(\varphi_{x}, \psi) \psi_{x}, 0, 0, 0, 0)^{t}  \nonumber \\
		&= (0, \sigma_{\varphi_{x}}(V_{x}^{1}, V^{3}) V_{xx}^{1} - kV_{xx}^{1} 
		+ \sigma_{\psi}(V_{x}^{1}, V^{3}) V_{x}^{3} - kV_{x}^{3}, 0, 0, 0, 0)^{t}.  \nonumber
	\end{align}
	
	Taking into account that $F(V, V_x)(t, \cdot) \in D(A)$ for $V \in \mathcal{H}_3$, 
	it follows from the Duhamel's principle that 
	\begin{align}
		V(t) = e^{tA} V_0 + \int^t_0 e^{(t-\tau)A} F(V, V_x)(\tau) \mathrm{d}r. \label{DUHAMMEL}
	\end{align}
	
	The existence of a global solution as well as its exponential decay can be
	proved as in \cite{MuRa2002} using a similar technique as for nonlinear Cauchy problems in \cite{Ra1992}.
	
	We assume that the initial data are small in the $\mathcal{H}_2$-norm, i.e. 
	\begin{align}
		\|V_0\|_{\mathcal{H}_2} < \delta.  \nonumber
	\end{align}
	
	Moreover, let us assume the boundness of $V_0$ in the $\mathcal{H}_3$-norm, i.e. let 
	\begin{align}
		\|V_0\|_{\mathcal{H}_3} < \nu  \nonumber
	\end{align}
	hold for a $\nu > 1$.
	
	Due to the smoothness of the solution, there exist two intervals $[0, T^0]$
	und $[0, T^1]$ such that 
	\begin{align}
		\|V(t)\|_{\mathcal{H}_2} &\leq \delta, \quad \forall t \in [0, T^0], \nonumber \\
		\|V(t)\|_{\mathcal{H}_3} &\leq \nu, \quad \forall t \in [0, T^1].  \nonumber
	\end{align}
	
	Let $d > 1$ be a constant to be fixed later on. We define two positive numbers $T^1_M$ und $T^0_M$ 
	as the biggest interval length such that the local solution satisfies 
	\begin{align}
		\|V(t)\|_{\mathcal{H}_2} \leq 2 c_1 \delta, \quad \forall t \in [0, T^0_M]  \nonumber
	\end{align}
	and 
	\begin{align}
		\|V(t)\|_{\mathcal{H}_3} \leq d\nu, \quad \forall t \in [0, T^1_M], \nonumber
	\end{align}
	respectively, fulfiling 
	\begin{align}
		\left\|e^{tA} V_0\right\|_{\mathcal{H}_2} \leq c_1 \|V\|_{\mathcal{H}_2} \nonumber
	\end{align}
	for the constant $c_1 > 0$ defined as in (\ref{ESTIMATE_HIGHER_ENERGIES}).
	
	Under these conditions, we obtain the following estimate for high energy.
	
	\begin{Lemma} \label{LEMMA_ENERGY_ESTIMATE} 
		There exist positive constants $c_{2},c_{3}$ independent of $V_{0}$ and $T$ such that the local solution 
		from Theorem \ref{LOCAL_NONLINEAR_EXISTENCE_THEOREM_HIGHER_REGULARITY} satisfies for $t\in[0, T_{M}^{1}]$ the inequality 
		\[
			\|V(t)\|_{\mathcal{H}_{3}}^{2}
			\leq c_{2}\|V_{0}\|_{\mathcal{H}^{3}}^{2}e^{c_{3}\sqrt{d\nu}\int_{0}^{t}\|V(\tau)\|_{\mathcal{H}_{2}}^{1/2}d\tau}.
		\]
	\end{Lemma}
	
	\begin{Proof}
		As our nonlinearity coincides with that considered for nonlinear Timoshenko
		systems with classical heat conduction and the estimates for the linear
		terms produced by our two new dissipations can be done in the same manner,
		we can repeat the proof from \cite{MuRa2002} literally.
	\end{Proof}
	
	Using Lemma \ref{LEMMA_ENERGY_ESTIMATE} and equality (\ref{DUHAMMEL}), we can write 
	\begin{align}
		F(V, V_x)(\tau) \in D(A) \subset \mathcal{H}_2, \quad \tau \geq 0,  \nonumber
	\end{align}
	we can estimate for $\|V(t)\|_{\mathcal{H}_2}$ 
	\begin{align}
		\|V(t)\|_{\mathcal{H}_2} &\leq \left\|e^{tA} V_0\right\|_{\mathcal{H}_2} 
		+ \int^t_0 \left\|e^{(t-\tau)A} F(V, V_x)(\tau)\right\|_{\mathcal{H}_2} \mathrm{d}\tau  \nonumber \\
		&\leq c_1 e^{-\alpha t} \|V_0\|_{\mathcal{H}_2} + c_1 \int^t_0 e^{-\alpha(t - \tau)} \left\|F(V, V_x)\right\|_{\mathcal{H}_2} \mathrm{d}\tau
		\label{LOCAL_SOLUTION_ESTIMATE}
	\end{align}
	by estimating the nonlinearity $F$ as in the lemma below.
	
	\begin{Lemma} \label{LEMMA_NONLINEARITY_ESTIMATE} 
		There exists a positive constant $c$
		such that the inequality 
		\[
			\|F(W,W_{x})\|_{\mathcal{H}_{2}}\leq c\|W\|_{\mathcal{H}_{2}}^{2}\|W\|_{\mathcal{H}_{3}}
		\]
		holds for all $W \in \mathcal{H}_{3}$ with $\|W\|_{\mathcal{H}_{2}} < C < \infty$.
	\end{Lemma}
	
	Further, we can show the following weighted a priori estimate.
	
	\begin{Lemma} \label{LEMMA_A_PRIORI_ESTIMATE} 
		Let 
		\[
			M_{2}(t) := \sup_{0 \leq \tau \leq t}\left(e^{\alpha \tau}\|V(\tau)\|_{\mathcal{H}_{2}}\right) 
		\]
		be defined for $t \in [0, T_{M}^{1}]$.
		
		There exist then $M_{0} > 0$ and $\delta > 0$ such that 
		\[
			M_{2}(t) \leq M_{0} < \infty 
		\]
		holds if $\|V_{0}\|_{\mathcal{H}_{3}} < \nu$ and $\|V_{0}\|_{\mathcal{H}_{2}} < \delta$,
		
		Moreover, $M$ does not depend on $T_{M}^{1}$ and $V_{0}$.
	\end{Lemma}
	
	\begin{Proof}
		We assume that $\|V\|_{\mathcal{H}_{2}}$ is bounded. Using Lemma \ref{LEMMA_NONLINEARITY_ESTIMATE} 
		and the estimate (\ref{LOCAL_SOLUTION_ESTIMATE}), we have 
		\[
			\|V(t)\|_{\mathcal{H}_{2}} \leq c_{1} e^{-\alpha t}\|V_{0}\|_{\mathcal{H}_{2}}
			+ c\int_{0}^{t} e^{-\alpha(t - \tau)}\|V(\tau )\|\| V(\tau )\|_{\mathcal{H}_{3}} \mathrm{d}\tau.
		\]
		
		With the aid of Lemma \ref{LEMMA_ENERGY_ESTIMATE}, there results 
		\begin{align}
			\|V(t)\|_{\mathcal{H}_{2}} \leq& c_{1}\|V_{0}\|_{\mathcal{H}_{2}}e^{-\alpha t}  \nonumber \\
			&+ c\int_{0}^{t} e^{-\alpha(t - \tau )}\|V(\tau )\|_{\mathcal{H}_{2}}^{2}\|V_{0}\|_{\mathcal{H}_{3}} 
			e^{c\sqrt{d\nu }\int_{0}^{\tau} \|V(r)\|_{\mathcal{H}_{2}}^{1/2} \mathrm{d}r} \mathrm{d}\tau. \label{V_H2_ESTIMATE}
		\end{align}
		
		Under assumption $\|V_{0}\|_{\mathcal{H}_{2}}\leq \delta$ for some $\delta >0$ to be determined later on, 
		we obtain for $t \in [0, \min\{T_{m}^{0}, T_{m}^{1}\}]$ 
		\begin{align}
			\|V(t)\|_{\mathcal{H}_{2}} &\leq c_{1} \delta e^{-\alpha t} 
			+ c\delta^{1/2} \nu e^{c\sqrt{d \nu} \int_{0}^{t }\|V(\tau )\|_{\mathcal{H}_{2}}^{1/2} \mathrm{d}\tau}
			\int_{0}^{t}e^{-\alpha(t - \tau)}\|V(\tau)\|_{\mathcal{H}_{2}}^{3/2} \mathrm{d}\tau   \nonumber \\
			&\leq c_{1} \delta e^{-\alpha t} + c\delta^{1/2} \nu e^{c\sqrt{d\nu}
			\int_{0}^{t} e^{-\alpha \tau /2} e^{\alpha \tau/2} \|V(\tau)\|_{\mathcal{H}_{2}}^{1/2} \mathrm{d} \tau} \times   \nonumber \\
			&{\phantom{=}} \times \int_{0}^{t}e^{-\alpha(t - \tau)}(e^{-\alpha \tau}
			e^{\alpha \tau }\|V(\tau )\|_{\mathcal{H}_{2}})^{3/2} \mathrm{d}\tau \nonumber \\
			&\leq c_{1} \delta e^{-\alpha t} + c\delta ^{1/2} \nu 
			e^{cte^{-\alpha t/2}\sqrt{d\nu }\sqrt{M(t)}}M_{2}(t)^{3/2}\int_{0}^{t} e^{-(\alpha - \tau)}
			e^{-3\alpha\tau/3} \mathrm{d}\tau  \nonumber \\
			&\leq c_{1} \delta e^{-\alpha t} + c \delta^{1/2} \nu e^{\sqrt{d\nu} \sqrt{M(t)}}M_{2}(t)^{3/2}\int_{0}^{t}
			e^{-(\alpha - \tau)}e^{-3\alpha \tau/3}\mathrm{d}\tau,  \nonumber
		\end{align}
		whence one can easily deduce 
		\[
			M_{2}(t)\leq c_{1}\delta + c\delta ^{1/2} \nu e^{c\sqrt{d\nu }\sqrt{M_{2}(t)}}M_{2}(t)^{3/2} \sup_{0\leq t < \infty}
			e^{\alpha t} \int_{0}^{t} e^{-\alpha(t - \tau)}e^{-3 \alpha \tau/2} \mathrm{d}\tau 
		\]
		after multiplication with $e^{\alpha t}$.
		
		From 
		\[
			\sup_{0\leq t < \infty} e^{\alpha t} \int_{0}^{t} e^{-\alpha (t-\tau)} e^{-3\alpha \tau/2}
			\mathrm{d} \tau = \sup_{0 \leq t < \infty} - \frac{2}{5} \frac{e^{-5\alpha \tau /2}}{\alpha} \big|_{\tau = 0}^{\tau = t}
			\leq c < \infty 
		\]
		it follows that 
		\[
			M_{2}(t) \leq c_{1} \delta + c \delta^{1/2} \nu M_{2}(t)^{3/2} e^{c\sqrt{d\nu}\sqrt{M_{2}(t)}}.
		\]
		
		We define a function 
		\[
			f(x) := c_{1} \delta + c \delta^{1/2} \nu x^{3/2} e^{c\sqrt{d\nu} \sqrt{x}} - x.
		\]
		
		We compute $f(0) = c_{1} \delta$ and $f^{\prime}(0) = -1$. According to the fundamental theorem of calculus, we know 
		\[
			f(x) = f(0) + \int_{0}^{x} f^{\prime}(\xi) \mathrm{d}\xi = c_{1} \delta + \int_{0}^{x} f^{\prime}(\xi)\mathrm{d}\xi.
		\]
		
		For sufficiently small $x$, we get $f^{\prime}(\xi) \leq - \frac{1}{2}$.
		This means that 
		\[
			f(x) \leq c_{1} \delta - \frac{1}{2}x.
		\]
		
		If we choose now a $\delta < \delta_{1} := \frac{x}{2c_{1}}$, we obtain $f(x) < 0$. 
		Because $f$ is continuous and $f(0) > 0$ as well as $f(x) < 0$ holds, $f$ must possess a zero in interval $[0, x]$. 
		Let $M_{0}$ be the smallest zero of $f$ in $[0, x]$. The latter must exist as $f^{-1}(\{0\}) \cap [0, x]$ is compact.
		
		We fix a $\delta_{2} < \delta_{1}$ to be so small that for $M_{2}(0) = \|V_{0}\|_{\mathcal{H}_{2}} < \delta_{2}$ 
		\[
			M_{2}(t) \leq M_{0}
		\]
		is fulfiled. It is possible due to the continuity of $M_{2}(t)$.
		
		Thus, $M_{2}(t)$ is bounded by a $M_{0}$ for all $t \in [0, \min\{T_{M}^{0}, T_{M}^{1}\}]$.
		
		If $T_{M}^{0} \geq T_{M}^{1}$, the claim of the theorem holds for $\delta < \delta_{2}$ und $M_{0}(\delta_{1}) < \infty $.
		
		Otherwise, we have $T_{M}^{0} < T_{M}^{1}$. We observe that for sufficiently small $\delta_{3} > 0$ 
		\[
			f(2c_{1}\delta ) = c\nu c_{1}^{3/2} e^{c\sqrt{d\nu} \sqrt{c_{1}\delta}}\delta^{2} - c_{1}\delta < 0
		\]
		is valid for $\delta < \delta_{3}$.
		
		We choose now an appropriately small $\delta_{3}$ to fulfil the above inequality.
		
		Hence, 
		\[
			\|V(t)\|_{\mathcal{H}_{3}} \leq M_{2}(t) \leq M_{0} < 2c_{1} \delta 
		\]
		for $\delta < \min\{\delta_{2}, \delta_{3}\}$.
		
		This contradicts to the maximality of $T_{M}^{0}$. That is why $T_{M}^{0} \geq T_{M}^{1}$ must be valid, 
		i.e. the claim holds for $\delta <\min\{\delta_{2}, \delta_{3}\}$ and $M_{0}(\delta_{1}) < \infty $.
	\end{Proof}
	
	This enables us finally to formulate and prove the theorem on global existence and exponential stability.
	
	\begin{Theorem}
		Let the assumptions of Theorem \ref{LOCAL_NONLINEAR_EXISTENCE_THEOREM_HIGHER_REGULARITY} be fulfiled. Moreover, let 
		\[
			\int_{0}^{L}\varphi_{0}(x)\mathrm{d}x = \int_{0}^{L}\varphi_{1}(x) \mathrm{d}x = \int_{0}^{L}\theta (x) \mathrm{d}x = 0.
		\]
		
		Let $\nu > 1$ be arbitrary but fixed. We can then find an appropriate $\delta > 0$ such that if $\|V_{0}\|_{\mathcal{H}_{2}} < \delta$ 
		and $\|V_{1}\|_{\mathcal{H}_{3}} < \nu$ hold there exists a unique global solution $(\varphi, \psi, \theta, q)$ 
		to (\ref{TIMOSHENKO_CATTANEO_NONLINEAR_DAMPED_EQUIV_1})---(\ref{TIMOSHENKO_CATTANEO_NONLINEAR_DAMPED_EQUIV_IC_1}) satisfying 
		\begin{align}
			\varphi, \psi &\in C^{3}([0, \infty )\times [0,L]),  \nonumber \\
			\theta, q &\in C^{2}([0, \infty )\times [0,L]).  \nonumber
		\end{align}
		
		There exists besides a constant $C_{0}(V_{0}) > 0$ such that for all $t \geq 0$ 
		\[
			\|V(t)\|_{\mathcal{H}_{2}} \leq C_{0}e^{-\alpha t}
		\]
		with $\alpha > 0$ from Theorem \ref{THEOREM_EXPONENTIAL_STABILITY_LINEAR_2} is valid.
	\end{Theorem}
	
	\begin{Proof}
		Theorem \ref{LOCAL_NONLINEAR_EXISTENCE_THEOREM_HIGHER_REGULARITY} guarantees
		the existence of a local solution with the regularity 
		\begin{align}
			\varphi,\psi &\in C^{3}([0, T] \times [0, L]),  \nonumber \\
			\theta,q &\in C^{2}([0, T] \times [0, L]).  \nonumber
		\end{align}
	
		Lemmata \ref{LEMMA_ENERGY_ESTIMATE} and \ref{LEMMA_A_PRIORI_ESTIMATE} suggest that 
		\begin{align}
			\|V(t)\|_{\mathcal{H}_{3}} &\leq c\|V_{0}\|_{\mathcal{H}_{3}}e^{\tilde{c}\sqrt{d\nu} 
			\int_{0}^{t}\|V(\tau )\|_{\mathcal{H}_{2}}^{1/2}\mathrm{d} \tau}  \nonumber \\
			&\leq c\|V_{0}\|_{\mathcal{H}_{3}}e^{\tilde{c} \sqrt{d \nu M_{0}}} 
			\leq ce\|V_{0}\|_{\mathcal{H}_{3}}, \quad t \leq T_{M}^{1} \leq T, \nonumber
		\end{align}
		where $\tilde{c} > 0$ and $\delta$ are chosen sufficiently small in order $\tilde{c} \sqrt{d \nu M_{0}} < 1$ is fulfiled.
		
		We put $d := ce$ and find 
		\[
			\|V(t)\|_{\mathcal{H}_{3}} \leq d\|V_{0}\|_{\mathcal{H}_{3}} < d \nu, \quad t \leq T_{M}^{1} \leq T.
		\]
		
		For $T_{M}^{1} < T$, we become a contradition to the maximality of $T_{M}^{1}$. Thus, $T_{M}^{1} = T$ must hold.
		
		If $0 \leq t \leq T$, there results from (\ref{V_H2_ESTIMATE}) that 
		\begin{align}
			\|V(t)\|_{\mathcal{H}_{2}} &\leq c_{1}\|V_{0}\|_{\mathcal{H}_{2}}e^{-\alpha t}
			+ c\int_{0}^{t}e^{-\alpha (t - \tau)}\|V(\tau )\|_{\mathcal{H}_{2}}^{2}\|V_{0}\|_{\mathcal{H}_{3}} e^{c\sqrt{d\nu }
			\int_{0}^{\tau}\|V(r)\|_{\mathcal{H}_{2}}^{1/2}\mathrm{d}r}\mathrm{d}\tau   \nonumber \\
			&\leq c\|V_{0}\|_{\mathcal{H}_{2}} + c\int_{0}^{t} \|V(\tau)\|_{\mathcal{H}_{2}}^{2}
			e^{c\sqrt{M_{2}(t)}} \mathrm{d}\tau   \nonumber \\
			&\leq c \|V_{0}\|_{\mathcal{H}_{2}} + ce^{c\sqrt{M_{2}(t)}}M_{2}(t) 
			\int_{0}^{t} \|V(\tau)\|_{\mathcal{H}_{2}} \mathrm{d}\tau, \nonumber
		\end{align}
		whence we conclude 
		\begin{equation}
			\|V(t)\|_{\mathcal{H}_{2}} \leq K \|V_{0}\|_{\mathcal{H}_{2}} \label{V_ESTIMATE_IN_H2}
		\end{equation}
		using the Gronwall's lemma for 
		\[
			K := cM_{0} e^{c (\sqrt{M_{0}} + M_{0})}.
		\]
		
		We choose $\delta^{\prime}$ such that $0 < \delta^{\prime} < \frac{\delta}{K}$ and obtain 
		\[
			\|V(T)\|_{\mathcal{H}_{2}} \leq K \|V_{0}\|_{\mathcal{H}_{2}} \leq K\delta^{\prime} \leq \delta.
		\]
		
		Therefore, there exists a continuation of $V$ onto $[T,T + T_{1}(\delta_{1})]$. 
		With (\ref{V_ESTIMATE_IN_H2}) there follows 
		\[
			\|V(T + T_{1}(\delta))\|_{\mathcal{H}_{2}} \leq K \|V_{0}\|_{\mathcal{H}_{2}} \leq \delta,
		\]
		i.e. we can smoothly continue the solution onto $[T+T_{1}(\delta_{1}), T + 2T_{1}(\delta_{1})]$.
		
		Here, we applied (\ref{V_ESTIMATE_IN_H2}) to the solution of the initial boundary value problem 
		with the initial value $W_{0} := V(T)$. This is allowed since $\|W_{0}\|_{\mathcal{H}_{3}} < \delta$ 
		and $\|W_{0}\|_{\mathcal{H}_{3}} \leq c < \infty $ holds according to Lemma \ref{LEMMA_ENERGY_ESTIMATE}.
		
		Hence, we can succesively obtain a global solution $V = (\varphi, \varphi_{t}, \psi, \psi_{t}, \theta, q)^{t}$ with 
		\begin{align}
			\varphi, \psi &\in C^{3}([0, \infty) \times [0, L]),  \nonumber \\
			\theta, q &\in C^{2}([0, \infty) \times [0, L]).  \nonumber
		\end{align}
		
		In particular, we can conlude 
		\[
			M_{2}(t) \leq M_{0} < \infty ,
		\]
		for all $t\in [0,\infty )$, since 
		\[
			\|V(t)\|_{\mathcal{H}_{2}} \leq K\delta^{\prime} \leq \delta .
		\]
		
		Finally, it follows 
		\[
			\|V(t)\|_{\mathcal{H}_{2}}\leq M_{0} e^{-\alpha t}.
		\]
	\end{Proof}

\vspace{0.3cm} \noindent \textbf{Acknownledgement} The author Salim A. Messaoudi has been funded during the work on this paper by KFUPM 
under Project \# SB070002.


\begin{thebibliography}{99}
\bibitem{Ti1921} Timoshenko S., {\it On the correction for shear of the differential equation for transverse vibrations of prismatic bars},
		Philosophical magazine {\bf 41} (1921), 744--746
		
	\bibitem{KiRe1987} Kim J.U. and Renardy Y., {\it Boundary control of the Timoshenko beam},
		SIAM J. Control Optim. {\bf 25} no. {\bf 6} (1987), 1417--1429
		
	\bibitem{FeShiZha1998} Feng D-X, Shi D-H, and Zhang W., {\it Boundary feedback stabilization of Timoshenko beam with boundary dissipation},
		Sci. China Ser. {\bf A 41} no. {\bf 5} (1998), 483--490
	
	\bibitem{RaFeSaCa2005} Raposo C.A., Ferreira J., Santos M.L., and Castro N.N.O.,
		{\it Exponential stability for Timoshenko system with two weak dampings},
		Applied Math Letters {\bf 18} (2005), 535--541
		
	\bibitem{Ta2000} Taylor S.W., {\it A smoothing property of a hyperbolic system and boundary controllability},
		J. Comput. Appl. Math. {\bf 114} (2000), 23--40
		
	\bibitem{LiuZhe1999} Liu Z. and Zheng S., {\it Semigroups associated with dissipative systems},
		Chapman \& Hall/CRC, 1999
		
	\bibitem{SouWeh2003} Soufyane A. and Wehbe A., {\it Uniform stabilization for the Timoshenko beam by a locally distributed damping},
		Electron. J. Differential Equations no. {\bf 29} (2003), 1--14
		
	\bibitem{Sou1999} Soufyane A., {\it Stabilisation de la poutre de Timoshenko},
		C.R. Acad. Sci. Paris S\'{e}r. I Math. {\bf 328} no. {\bf 8} (1999), 731--734
		
	\bibitem{ShiFe2002} Shi D-H and Feng D-X, {\it Exponential decay rate of the energy of a Timoshenko beam with locally distributed feedback},
		ANZIAM J. {\bf 44} no. {\bf 2} (2002), 205--220
		
	\bibitem{MuRa2007} Mu\~{n}oz Rivera J.E. and Racke R., {\it Timoshenko systems with indefinite damping},
		J. Math. Anal. Appl. {\bf 341} (2008), 1068--1083
		
	\bibitem{MuRa2003} Mu\~{n}oz Rivera J.E. and Racke R., {\it Global stability for damped Timosheko systems},
		Discrete Contin. Dyn. Syst. {\bf 9} no. {\bf 6} (2003), 1625--1639
		
	\bibitem{XuYu2003} Xu G-Q and Yung S-P, {\it Stabilization of Timoshenko beam by means of pointweise controls}
		ESAIM, Control Optim. Calc. Var. {\bf 9} (2003), 579--600
	
	\bibitem{AmBeRiRa2003} Ammar-Khodja F., Benabdallah A., Mu\~{n}oz Rivera J.E. and Racke R., 
		{\it Energy decay for Timoshenko systems of memory type},
		J. Differential Equations {\bf 194} no. {\bf 1} (2003), 82--115
		
	\bibitem{Sa2002} Santos M.L., {\it Decay rates for solutions of a Timoshenko system with a memory condition at the boundary},
		Abstr. Appl. Anal. {\bf 7} no. {\bf 10} (2002), 531--546
		
	\bibitem{ShiFe2001} Shi D-H and Feng D-X, {\it Exponential decay of Timoshenko beam with locally distributed feedback},
		IMA J. Math. Control Inform. {\bf 18} no. {\bf 3} (2001), 395--403
	
	\bibitem{MuRa2002} Mu\~{n}oz Rivera J.E. and Racke R., 
		{\it Mildly dissipative nonlinear Timoshenko systems --- global existence and exponential stability},
		J. Math. Anal. Appl. {\bf 276} (2002), 248--278
		
	\bibitem{Ta1992} Tarabek M.A., {\it On the existence of smooth solutions in one-dimensional thermoelasticity with second sound},
		Quart. Appl. Math. {\bf 50} (1992), 727--742
		
	\bibitem{Ra2002} Racke R., {\it Thermoelasticity with second sound --- exponential stability in linear and nonlinear 1-d},
		Math. Meth. Appl. Sci. {\bf 25} (2002), 409--441
		
	\bibitem{MeSa2005} Messaoudi S.A. and Said-Hourari B., {\it Exponential stability in one-dimensional nonlinear thermoelasticity with second sound},
		Math Methods Appl. Sci. {\bf 28} (2005), 205--232
		
	\bibitem{Ra2003} Racke R., {\it Asymptotic behaviour of solutions in linear 2- or 3-d thermoelasticity with second sound},
		Quarterly of Applied Math. {\bf 61} \# {\bf 2} (2003), 315--328
		
	\bibitem{Me2002} Messaoudi S.A., {\it Local existence and blow up in thermoelasticity with second sound},
		CPDE Vol. {\bf 26} \# {\bf 8} (2002), 1681--1693
		
	\bibitem{MeSa2004} Messaoudi S.A. and Said-Hourari B., 
		{\it Blow up of solutions with positive energy in nonlinear thermo-elasticity with second sound},
		J. Appl. Math. {2004} \# {\bf 3} (2004), 201--211

	\bibitem{Ir2006} Irmscher T., {\it Aspekte hyperbolischer Thermoelastizit\"{a}t},
		Dissertation, Konstanz (2006)
			
	\bibitem{FeRa2007} Fern\'{a}ndez Sare H.D. and Racke R., {\it On the stability of damped Timoshenko Systems --- Cattaneo versus Fourier's law},
		Konstanzer Schriften Math. Inf. {\bf 227} (2007)
		
	\bibitem{Sl1981} Slemrod M., {\it Global existence, uniqueness, and asymptotic stability
		of classical smooth solutions in one-dimensional non-linear thermoelasticity}, Arch. Rational Mech. Anal. {\bf 76} (1981), 97-133
		
	\bibitem{JiRa2000} Jiang S., Racke R., {\it Evolution Equations in Thermoelasticity},
    		Chapman \& Hall/CRC, Monographs and Surveys in Pure and Applied Mathematics, Vol. 112 (2000)
		
	\bibitem{Ra1992} Racke R., {\it Lectures on Nonlinear Evolution Equations. Initial Value Problems},
    		Aspects of Math., Vol. E19, Friedr. Vieweg and Sohn, Braunschweig/Wiesbaden (1992)
\end{thebibliography}
\end{document}